\documentclass[11pt,a4paper]{article}
\usepackage[utf8]{inputenc}
\usepackage[T1]{fontenc}
\usepackage{amsfonts}
\usepackage{amssymb}
\usepackage{amsthm}
\usepackage{amsmath}
\usepackage{graphicx}
\usepackage{mathrsfs}
\usepackage[left=2cm,right=2cm,top=2cm,bottom=2cm]{geometry}

\newtheorem{theorem}{Theorem}
\newtheorem{definition}[theorem]{Definition}

\newtheorem{lemma}[theorem]{Lemma}

\newtheorem{proposition}[theorem]{Proposition}
\newtheorem{remark}[theorem]{Remark}

\begin{document}

\title{ Mean Field Game for Linear Quadratic Stochastic Recursive Systems }
\author{Liangquan Zhang$^{1}$\thanks{
L. Zhang acknowledges the financial support partly by the National Nature
Science Foundation of China (Grant No. 11701040, 61871058, 11871010 \&
61603049) and the Fundamental Research Funds for the Central Universities
(No. 500417024 \& 505018304). E-mail: xiaoquan51011@163.com},\ Xun Li$^{2}$%
\thanks{%
X. Li acknowledges the financial support partly by PolyU G-UA4N, Hong Kong
RGC under grants 15224215 and 15255416. E-mail: li.xun@polyu.edu.hk.} \  \\
{\small 1. School of Science }\\
{\small \ \ Beijing University of Posts and Telecommunications }\\
{\small \ \ Beijing 100876, China }\\
{\small 2. Department of Applied Mathematics }\\
{\small \ The Hong Kong Polytechnic University, Hong Kong}\\
}
\maketitle

\begin{abstract}
This paper focuses on linear-quadratic (LQ for short) mean-field games
described by forward-backward stochastic differential equations (FBSDEs for
short), in which the individual control region is postulated to be convex.
The decentralized strategies and consistency condition are represented by a
kind of coupled mean-field FBSDEs with projection operators. The
well-posedness of consistency condition system is obtained using the
monotonicity condition method. The $\epsilon$-Nash equilibrium property is
discussed as well.
\end{abstract}



\section{Introduction}

\label{sect:1}

The control of stochastic multi-agent systems has attracted large attentions
by many researchers. As well-known, the large population systems arise
naturally in various different fields (e.g., biology, engineering, social
science, economics and finance, operational research and management, etc.).
Readers interested in this topic may refer \cite{HM, HCM, HCM2, HMC} for
more details of their solid backgrounds and real applications. The agents
(or players) in large population system are individually negligible but
their collective behaviors will make some significant impact on all agents.
This trait can be captured by the weakly-coupling structure in the
individual dynamics and cost functionals through the state-average. The
individual behaviors of all agents in micro-scale can make their mass
effects in the macro-scale.

As for the controlled large population system, it is intractable for a given
agent to collect all agents due to the highly complex interactions among its
colleagues. Consequently, the centralized controls, which are built upon the
full information of all agents, are not implementable and not efficient in
large population framework. Alternatively, it is more reasonable and
effective to study the decentralized strategies which depend on the local
information\footnote{%
Here local information means the optimal control regulator for a given
agent, is designed on its own individual state and some quantity which can
be obtained in off-line way.} only. The mean-field type stochastic control
problem is of both great interest and importance in various fields such as
science, engineering, economics, management, and particularly in financial
investment. In contrast with the standard stochastic control problems, the
underlying dynamic system and the cost functional involve state processes as
well as their expected values (hence the name mean-field). In financial
investment, however, one frequently encounters interesting problems which
are closely related to money managers' performance evaluation and incentive
compensation mechanisms. Together with MF-FSDEs, research is naturally
required on optimal control problems based on mean-field forward-backward
stochastic differential equations (hereafter MF-FBSDEs). Hence, one powerful
tool employed is so-called mean-field games (see \cite{LL}). The basic idea
is to approximate the initial large population control problem by its
limiting problem via some mean-field term (i.e., the asymptotic limit of
state-average). There are huge literature can be found in \cite{Bardi, CD,
cfs, pc, GLL, KTY, LZ, LL} for the study of mean-field games; \cite{HCM2}
for cooperative social optimization; \cite{HM}, \cite{NH} and \cite{NC} and
references therein for models with a major player; \cite{AD, BDL, Y} for
optimal control with a mean term in the dynamics and cost, etc.

The main contribution of this paper is to study the forward backward
mean-field LQG of large population systems for which the individual states
follow some forward backward stochastic differential equations (FBSDEs in
short). This framework makes our setting very different to existing works of
mean-field LQG games wherein the individual states evolve by some forward
stochastic differential equations. In contrast to classical stochastic
differential equations, the terminal condition of BSDE should be specified
as the priori random variable, which means, the BSDE will admit one pair of
adapted solutions, in which the second solution component (the diffusion
term) is naturally appeared here by virtue of the martingale representation
theorem and the adaptiveness requirement for filtration. The linear BSDEs
are first introduced in \cite{Bis} for studying the optimal control
problems, and the general nonlinear BSDEs are developed by Pardoux and Peng
in 1992 \cite{PP}. Since then, the study of BSDE has initiated consistent
and intense discussions, moreover, it has been used in many applications of
diverse areas. For instance, the BSDE takes very important role to
characterize the nonlinear expectation ($g$-expectation, see \cite{Pg}), or
the stochastic differential recursive utility (see \cite{DE}). Subsequently,
El Karoui, Peng, and Quenez \cite{KPQ} presents many applications of BSDE in
mathematical finance and optimal control theory. Pardoux and Peng establish
a kind of stochastic partial differential equations with backward doubly SDE
(see \cite{PP1994}). Therefore, it is very natural to study its dynamic
optimization in large-population setting. Indeed, the dynamic optimization
of backward large population system is inspired by a variety of scenarios.
For example, the dynamic economic models for which the participants are of
some recursive utilities or nonlinear expectations, or some production
planning problems with some tracking terminal objectives but affected by the
market price via production average.

Another example arises from the risk management when considering the
relative or comparable criteria based on the average performance of all
other peers through the whole sector. This is the case for a given pension
fund to evaluate its own performance by setting the average performance
(say, average hedging cost or initial deposit, surplus) as its benchmark. In
addition, the controlled forward large population systems, which are
subjected to some terminal constraints, can be reformulated by some backward
large population systems, as motivated by \cite{LimZ}. Applying to
performance evaluation and incentive compensation of fund managers in the
field of financial engineering is of both academic and practical importance.
Findings from FBSDEs and mean-field stochastic optimal controls will not
only contribute to the academic literature by shedding light on performance
evaluation and incentive compensation schemes, but also provide practical
applications to fund management and risk control, especially under the
current circumstances with on-going economic recession. More importantly,
research outcomes in this field are expected to add to our knowledge based
on economic theory about providing appropriate incentives for managers in an
agency framework. They can also be generalized to various industries and
economic regions to provide policy makers with a theoretical basis during
their decision-making processes. Inspired by above mentioned motivations,
this paper studies the forward backward mean-field linear-quadratic-Gaussian
(BMFLQG) games.

We concern on the linear-quadratic (LQ) mean-field game where the individual
control domain is convex subset of $\mathbb{R}^{m}$. The LQ problems with
convex control domain comes naturally from various practical applications.
For instance, the no-shorting constraint in portfolio selection leads to the
LQ control with positive control ($\mathbb{R}_{+}^{m},$ the positive
orthant). Moreover, due to general market accessibility constraint, it is
also interesting to study the LQ control with more general closed convex
cone constraint (see \cite{hz}). As a response, this paper investigates the
LQ dynamic game of large-population system with general closed convex
control constraint.

The control constraint brings some new features to our study here: (1) The
related consistency condition (CC) system is no longer linear, and it
becomes a class of nonlinear FBSDEs with projection operator. (2) Due to the
nonlinearity of (1), the standard Riccati equation with feedback control is
no longer valid to represent the consistency condition of limit
state-average process. Instead, the consistency condition is embedded into a
class of mean-field coupled FBSDEs with a generic driven Brownian motion.

Similarly like in Hu, Huang, Li \cite{HHL}, we first apply the stochastic
maximum principle for convex control domain of the optimal decentralized
response through some Hamiltonian system with projection operator upon the
control set $U$ for forward-backward systems. Then, the consistency
condition system is connected to the well-posedness of some mean-field
coupled forward-backward stochastic differential equation (MF-FBSDE). Next,
we state some monotonicity condition of this MF-FBSDE to obtain its
uniqueness and existence. At last, the related approximate Nash equilibrium
property is also established. The MFG strategy derived is an open-loop
manner. Consequently, the approximate Nash equilibrium property is verified
under the open-loop strategies perturbation and some estimates of
forward-backward SDE are involved. In addition, all agents are set to be
statistically identical thus the limiting control problem and fixed-point
arguments are given for a representative agent.

In order to make our paper more accessible to the reader, we provide the
standard procedure of MFG, and describe our result mainly consisting of the
following steps:

\textit{Step 1}: Fix the state-average limit: lim$_{N\rightarrow +\infty
}x^{(N)}$ and by lim$_{N\rightarrow +\infty }y^{(N)}$ a frozen process $%
\mathbb{E}x$ and $\mathbb{E}y$ (the law of large numbers) and formulate an
auxiliary stochastic control problem for $\mathcal{A}_{i}$ which is
parameterized by $\mathbb{E}x$ and $\mathbb{E}y$. Note that the coupled
Hamiltonian systems admits a unique strong adapted solution.

\textit{Step 2:} Solve the above auxiliary stochastic control problem via
Pontryagin's maximum principle to obtain the decentralized optimal state $%
\left( x_{t}^{i,\ast },y_{t}^{i,\ast },z_{t}^{i,\ast }\right) $, (which
should depend on the undetermined processes $\mathbb{E}x$ and $\mathbb{E}y$%
). By means of convex analysis, we are able to construct the \textit{unique}
feedback control using a so-called projection mapping, $\varphi \left(
t,p,q,k\right) .$

\textit{Step 3:} Characterize the decentralized strategies $\{\bar{u}%
_{t}^{i},1\leq i\leq N\}$ of Problem (\textbf{CC}) through the auxiliary (%
\textbf{LCC}) and consistency condition system. Namely, show $\{\bar{u}%
_{t}^{i},1\leq i\leq N\}$ is an $\epsilon $-Nash equilibrium. This step
actually can be further divided into:

\textit{Step 3-1: }Introduce the decentralized state $\left( \breve{x}%
_{t}^{i},\breve{y}^{i},\breve{z}^{i}\right) ,$ with its decentralized
open-loop optimal strategy $\varphi \left( \chi ^{i},\beta ^{i},\gamma
^{i}\right) $ and the consistency conditions systems. We get two estimations
between them in Lemma \ref{lemma for xN and m} and Lemma \ref{lemma for x
and xi}, respectively;

\textit{Step 3-2: }For any fixed $i$, $1\leq i\leq N$, we shall consider a
group of state equations $\left( \tilde{x}^{i},\tilde{y}^{i},\tilde{z}%
^{i}\right) $ driven by certain perturbation control $u^{i}\in \mathcal{U}%
_{ad}^{d,i}$ and systems $\left( \mathring{x}^{i},\mathring{y}^{i},\mathring{%
z}^{i}\right) $ of the decentralized limiting state with perturbation
control. Similarly, we have the estimations between perturbation systems and
consistency condition system in Lemma \ref{lemma 2}, plus the estimations
between perturbation systems and decentralized limiting state with
perturbation control system in Lemma \ref{lemma 3};

\textit{Step 3-3: }Finally, based on Lemma \ref{lemma for xN and m}-Lemma %
\ref{lemma 3}, employing the relation between limiting cost functional $%
J_{i} $ and the cost functional $\mathcal{J}_{i}$ of $\mathcal{A}_{i}$ with
help of perturbational control, we are able to get our desired result.

This paper is organized as follows: Section \ref{sect2} formulates the LQ
MFGs of BSDEs type with convex control domain. The decentralized strategies
are derived with the help of a mean field forward-backward SDEs with
projection operators. The consistency condition is also established. Section %
\ref{sect3} verifies the $\epsilon $-Nash equilibrium of the decentralized
strategies. Some proofs will be scheduled in Appendix \ref{app1}. Related
results on properties of projection in convex analysis can be found in
Appendix \ref{app2}.

\section{Preliminaries}

\label{sect2}

Throughout this paper, we denote the $k$-dimensional Euclidean space by $%
\mathbb{R}^{k}$ with standard Euclidean norm $|\cdot |$ and standard
Euclidean inner product $\langle \cdot ,\cdot \rangle $. The transpose of a
vector (or matrix) $x$ is denoted by $x^{T}$. $\text{Tr}(A)$ denotes the
trace of a square matrix $A$. Let $\mathbb{R}^{m\times n}$ be the Hilbert
space consisting of all ($m\times n$)-matrices with the inner product $%
\langle A,B\rangle :=\text{Tr}(AB^{\top })$ and the norm $|A|:=\langle
A,A\rangle ^{\frac{1}{2}}$. Denote the set of symmetric $k\times k$ matrices
with real elements by $S^{k}$. If $M\in S^{k}$ is positive (semi)definite,
we write $M>\ (\geq )\ 0$. $L^{\infty }(0,T;\mathbb{R}^{k})$ is the space of
uniformly bounded $\mathbb{R}^{k}-$valued functions. If $M(\cdot )\in
L^{\infty }(0,T;S^{k})$ and $M(t)>\ (\geq )\ 0$ for all $t\in \lbrack 0,T]$,
we say that $M(\cdot )$ is positive (semi) definite, which is denoted by $%
M(\cdot )>\ (\geq )\ 0$. $L^{2}(0,T;\mathbb{R}^{k})$ is the space of all $%
\mathbb{R}^{k}-$valued functions satisfying $\int_{0}^{T}|x(t)|^{2}dt<\infty
.$

Consider a finite time horizon $[0,T]$ for fixed $T>0$. We assume $(\Omega ,%
\mathcal{F},\{\mathcal{F}_{t}\}_{0\leq t\leq T},P)$ is a complete, filtered
probability space on which a standard $N$-dimensional Brownian motion $%
\{W_{i}(t),\ 1\leq i\leq N\}_{0\leq t\leq T}$ is defined. For given
filtration $\mathbb{F}=\{\mathcal{F}_{t}\}_{0\leq t\leq T},$ let $L_{\mathbb{%
F}}^{2}(0,T;\mathbb{R}^{k})$ denote the space of all $\mathcal{F}_{t}$%
-progressively measurable $\mathbb{R}^{k}$-valued processes satisfying $%
\mathbb{E}\int_{0}^{T}|x(t)|^{2}dt<\infty .$ Let $L_{\mathbb{F}}^{2,\mathcal{%
E}_{0}}(0,T;\mathbb{R}^{k})\subset L_{\mathbb{F}}^{2}(0,T;\mathbb{R}^{k})$
be the subspace satisfying $\mathbb{E}x_{t}\equiv 0$ for $x\in L_{\mathbb{F}%
}^{2,\mathcal{E}_{0}}(0,T;\mathbb{R}^{k}).$ Let $L_{\mathbb{F}_{T}}^{2}(%
\mathbb{R}^{k})$ denote the space of all $\mathcal{F}_{T}$-measurable $%
\mathbb{R}^{k}$-valued random variable satisfying $\mathbb{E}|\xi
|^{2}<\infty .$

Now let us consider a large-population system with $N$ weakly-coupled
negligible agents $\{\mathcal{A}_{i}\}_{1\leq i\leq N}$. The state $x^{i}$
and $y^{i}$ for each $\mathcal{A}_{i}$ satisfies the following controlled
linear stochastic system:
\begin{equation}
\left\{
\begin{array}{lll}
\mathrm{d}x_{t}^{i} & = & \left(
A_{t}x_{t}^{i}+B_{t}u_{t}^{i}+F_{t}x_{t}^{(N)}+b_{t}\right) \mathrm{d}%
t+\left( D_{t}u_{t}^{i}+\sigma _{t}\right) \mathrm{d}W_{t}^{i}, \\
\mathrm{d}y_{t}^{i} & = & -\left(
M_{t}x_{t}^{i}+U_{t}y_{t}^{i}+H_{t}x_{t}^{(N)}+V_{t}y_{t}^{(N)}+K_{t}u_{t}^{i}+f_{t}\right)
\mathrm{d}t+z_{t}^{i}\mathrm{d}W_{t}^{i}, \\
x_{0}^{i} & = & x\in \mathbb{R}^{n},\text{ }y_{T}^{i}=\Phi x_{T}^{i},\qquad
0\leq t\leq T,%
\end{array}%
\right.  \label{FB1}
\end{equation}%
where $x^{(N)}(\cdot )=\displaystyle\frac{1}{N}\sum_{i=1}^{N}x^{i}(\cdot )$
and $y^{(N)}(\cdot )=\displaystyle\frac{1}{N}\sum_{i=1}^{N}y^{i}(\cdot )$ is
the state-average, $(A,B,F,b,D,\sigma,$\newline
$M,N,H,V,U,K,f,\Phi )$ are matrix-valued functions with appropriate
dimensions to be identified. For sake of presentation, we set all agents are
homogeneous or statistically symmetric with same coefficients $%
(A,B,F,b,D,\sigma ,M,U,H,V,U,K,f,\Phi )$ and deterministic initial states $x$%
.


Now we identify the information structure of large population system: $%
\mathbb{F}^{i}=\{\mathcal{F}^{i}_t\}_{0 \leq t \leq T}$ is the natural
filtration generated by $\{W_i(t), 0 \leq t \leq T\}$ and augmented by all $%
P-$null sets in $\mathcal{F}.$ $\mathbb{F}=\{\mathcal{F}_t\}_{0 \leq t \leq
T}$ is the natural filtration generated by $\{W_i(t), 1 \leq i \leq N, 0
\leq t \leq T\}$ and augmented by all $P-$null sets in $\mathcal{F}.$ Thus, $%
\mathbb{F}^{i}$ is the individual decentralized information of $i^{th}$
Brownian motion while $\mathbb{F}$ is the centralized information driven by
all Brownian motion components. Note that the heterogeneous noise $W_i$ is
specific for individual agent $\mathcal{A}_{i}$ but $x^i(t)$ is adapted to $%
\mathcal{F}_t$ instead of $\mathcal{F}^{i}_t$ due to the coupling
state-average $x^{(N)}.$

The (centralized) admissible control $u^{i}\in \mathcal{U}_{ad}^{c}$ where
the (centralized) admissible control set $\mathcal{U}_{ad}^{c}$ is defined
as
\begin{equation*}
\mathcal{U}_{ad}^{c}:=\{u^{i}(\cdot )|u^{i}(\cdot )\in L_{\mathbb{F}%
}^{2}(0,T;U),\quad 1\leq i\leq N\},
\end{equation*}%
where $U\subset \mathbb{R}^{m}$ is a closed convex set. Typical examples of
such set is $U=\mathbb{R}_{+}^{m}$ which represents the positive control. By
\textquotedblleft centralized\textquotedblright , we mean $\mathbb{F}$ is
the centralized information generated by all Brownian motion components.
Moreover, we also define decentralized control as $u^{i}\in \mathcal{U}%
_{ad}^{d,i}$, where the decentralized admissible control set $\mathcal{U}%
_{ad}^{d,i}$ is defined as
\begin{equation*}
\mathcal{U}_{ad}^{d,i}:=\{u^{i}(\cdot )|u^{i}(\cdot )\in L_{\mathbb{F}%
^{i}}^{2}(0,T;U),\quad 1\leq i\leq N\}.
\end{equation*}%
Note that both $\mathcal{U}_{ad}^{d,i}$ and $\mathcal{U}_{ad}^{c}$ are
defined in open-loop sense, and $\mathcal{U}_{ad}^{d,i}\subset \mathcal{U}%
_{ad}^{c}$. Let $u=(u^{1},\cdots ,u^{i}$, $\cdots ,u^{N})$ denote the set of
control strategies of all $N$ agents and $u_{-i}=(u^{1},\cdots ,u^{i-1},$ $%
u^{i+1},\cdots ,u^{N})$ denote the control strategies set except the $i^{th}$
agent $\mathcal{A}_{i}.$ Introduce the cost functional of $\mathcal{A}_{i}$
as
\begin{eqnarray}
\mathcal{J}_{i}\left( u^{i},u_{-i}\right) &=&\frac{1}{2}\mathbb{E}\Bigg [%
\int_{0}^{T}\Big [\left\langle Q_{t}\left( x_{t}^{i}-x_{t}^{(N)}\right)
,x_{t}^{i}-x_{t}^{(N)}\right\rangle +L_{t}\left(
y_{t}^{i}-y_{t}^{(N)}\right) ^{2}  \notag \\
&&+\left\langle R_{t}u_{t}^{i},u_{t}^{i}\right\rangle \Big ]\mathrm{d}%
t+\left\langle G\left( x_{T}^{i}-x_{T}^{(N)}\right)
,x_{T}^{i}-x_{T}^{(N)}\right\rangle \Bigg ].  \label{cost1}
\end{eqnarray}%
We assume the followng conditions:

\begin{description}
\item[(A1)] Assume $A(\cdot ),$ $F(\cdot ),$ $M(\cdot ),$ $U(\cdot ),$ $%
H(\cdot ),$ $V(\cdot )\in L^{\infty }(0,T;S^{n}),$ and $B(\cdot ),$ $D(\cdot
),$ $K\left( \cdot \right) \in L^{\infty }(0,T;\mathbb{R}^{n\times m}),$ $%
b(\cdot ),$ $\sigma (\cdot ),$ $f\left( \cdot \right) \in L^{\infty }(0,T;%
\mathbb{R}^{n});$

\item[(A2)] $Q(\cdot ),$ $L\left( \cdot \right) \in L^{\infty }(0,T;S^{n}),$
$Q(\cdot ),$ $L\left( \cdot \right) \geq 0,$ $R(\cdot )\in L^{\infty
}(0,T;S^{m}),$ $R(\cdot )>0$ and $R^{-1}(\cdot )\in L^{\infty }(0,T;S^{m})$,
$G\in S^{n}$, $G>0$.
\end{description}

By virtue of the theory of mean field BSDEs (see Lemma 3.1 in \cite{BDLP}),
under the assumptions (A1)-(A2), Eq. $(\ref{FB1})$ admits a unique solution $%
\left( x^{i},y^{i},z^{i}\right) \in L_{\mathbb{F}^{W}}^{2}(0,T;\mathbb{R}%
^{n})$ $\times L_{\mathbb{F}^{W}}^{2}(0,T;\mathbb{R}^{n})\times L_{\mathbb{F}%
^{W}}^{2}(0,T;\mathbb{R}^{n})$ with an admissible control $u_{i}\in \mathcal{%
U}_{ad}^{c}$. We now formulate the large population LQG with control
constraint (\textbf{CC}).

\textbf{Problem (CC).} Find an open-loop Nash equilibrium strategies set $%
\bar{u}=(\bar{u}^{1},\bar{u}^{2},\cdots ,\bar{u}^{N})$ satisfying
\begin{equation*}
\mathcal{J}_{i}(\bar{u}^{i}(\cdot ),\bar{u}_{-i}(\cdot ))=\inf_{u^{i}(\cdot
)\in \mathcal{U}_{ad}^{c}}\mathcal{J}_{i}(u^{i}(\cdot ),\bar{u}_{-i}(\cdot
)),
\end{equation*}%
where $\bar{u}_{-i}$ represents $(\bar{u}^{1},\cdots ,\bar{u}^{i-1},\bar{u}%
^{i+1},\cdots ,\bar{u}^{N})$, the strategies of all agents except $\mathcal{A%
}_{i}$.

Observe that the problem (\textbf{CC}) is of large computational issue since
the highly-complicated coupling structure among these agents. Alternatively,
the mean-field game theory employed is to search the approximate Nash
equilibrium, which bridges the \textquotedblleft centralized" LQG to the
limiting LQG control problems, as the number of agents tends to infinity.
Similar in \cite{HHL}, we need to construct some auxiliary control problem
using the frozen state-average limit. Based on it, we can find the
decentralized strategies by consistency condition.

Let us introduce the following auxiliary problem for $\mathcal{A}_{i}:$
\begin{equation}
\left\{
\begin{array}{lll}
\mathrm{d}x_{t}^{i,\diamond } & = & \left( A_{t}x_{t}^{i,\diamond
}+B_{t}u_{t}^{i}+F_{t}\phi _{t}^{1}+b_{t}\right) \mathrm{d}t+\left(
D_{t}u_{t}^{i}+\sigma _{t}\right) \mathrm{d}W_{t}^{i}, \\
\mathrm{d}y_{t}^{i,\diamond } & = & -\left( M_{t}x_{t}^{i,\diamond
}+U_{t}y_{t}^{i,\diamond }+H_{t}\phi _{t}^{1}+V_{t}\phi
_{t}^{2}+K_{t}u_{t}^{i}+f_{t}\right) \mathrm{d}t+z_{t}^{i,\diamond }\mathrm{d%
}W_{t}^{i}, \\
x_{0}^{i,\diamond } & = & x\in \mathbb{R}^{n},\text{ }y_{T}^{i,\diamond
}=\Phi x_{T}^{i,\diamond },\qquad 0\leq t\leq T,%
\end{array}%
\right.  \label{lcccontrol}
\end{equation}%
with the limiting cost functional given by
\begin{eqnarray}
J_{i}\left( u^{i}\right) &=&\frac{1}{2}\mathbb{E}\Bigg [\int_{0}^{T}\Big [%
\left\langle Q_{t}\left( x_{t}^{i,\diamond }-\phi _{t}^{1}\right)
,x_{t}^{i,\diamond }-\phi _{t}^{1}\right\rangle +L_{t}\left(
y_{t}^{i,\diamond }-\phi _{t}^{2}\right) ^{2}  \notag \\
&&+\left\langle R_{t}u_{t}^{i},u_{t}^{i}\right\rangle \Big ]\mathrm{d}%
t+\left\langle G\left( x_{T}^{i,\diamond }-\phi _{T}^{1}\right)
,x_{T}^{i,\diamond }-\phi _{T}^{1}\right\rangle \Bigg ],  \label{cost2}
\end{eqnarray}%
where $\phi ^{i},$ $i=1,2$ are the average limit of realized states which
should be determined by the consistency-condition (CC) in our later analysis
(see \eqref{cc}). Note that the auxiliary state $\left( x_{t}^{i,\diamond
},y_{t}^{i,\diamond },z_{t}^{i,\diamond }\right) $ is different to the true
state $\left( x^{i},y^{i},z^{i}\right) $. Also, the admissible control $%
u^{i} $ in \eqref{lcccontrol}, \eqref{cost2} $\in \mathcal{U}_{ad}^{d,i}$
whereas in \eqref{FB1}, \eqref{cost1}, the admissible control $\in \mathcal{U%
}_{ad}^{c}$ (for sake of simplicity, we still denote them with the same
notation).

Now we formulate the following limiting stochastic optimal control problem
with control constraint (\textbf{LCC}).

\textbf{Problem (LCC).} For the $i^{th}$ agent, $i=1,2,\cdots ,N,$ find $%
u^{i,\ast }(\cdot )\in \mathcal{U}_{ad}^{d,i}$ satisfying
\begin{equation*}
J_{i}(u^{i,\ast }(\cdot ))=\inf_{u^{i}(\cdot )\in \mathcal{U}%
_{ad}^{d,i}}J_{i}(u^{i}(\cdot )).
\end{equation*}%
Then $u^{i,\ast }(\cdot )$ is called a decentralized optimal control for
Problem (\textbf{LCC}). Now we apply the well known maximum principle
(Theorem 3.3 in \cite{Wu}) to characterize $u^{i,\ast }$ with the optimal
state ${x}^{i,\ast }.$ To this end, let us introduce the following adjoint
process%
\begin{equation*}
\left\{
\begin{array}{rcl}
\mathrm{d}p_{t}^{i} & = & \left[ U_{t}p_{t}^{i}-L_{t}\left( y_{t}^{i,\ast
}-\phi _{t}^{2}\right) \right] \mathrm{d}t, \\
\mathrm{d}q_{t}^{i} & = & \left[ -M_{t}p_{t}^{i}-A_{t}q_{t}^{i}+Q_{t}\left(
x_{t}^{i,\ast }-\phi _{t}^{1}\right) \right] \mathrm{d}t+k_{t}^{i}\mathrm{d}%
W_{t}^{i}, \\
p_{0}^{i} & = & 0,\text{ }q_{T}^{i}=\Phi ^{T}p_{T}^{i}-G\left( x_{T}^{i,\ast
}-\phi _{T}^{1}\right) .%
\end{array}%
\right.
\end{equation*}%
The Hamiltonian function can be expressed by
\begin{eqnarray}
\mathcal{H}^{i} &\triangleq &\mathcal{H}^{i}\left(
t,p^{i},q^{i},k^{i},x^{i},y^{i},z^{i},u^{i}\right)  \notag \\
&=&\left\langle p^{i},Mx^{i}+Uy^{i}+H\phi ^{1}+V\phi
^{2}+Ku^{i}+f\right\rangle  \notag \\
&&+\left\langle q^{i},Ax^{i}+Bu^{i}+F\phi ^{1}+b\right\rangle +\left\langle
k^{i},Du^{i}+\sigma \right\rangle  \notag \\
&&-\frac{1}{2}\left\langle Q\left( x^{i}-\phi ^{1}\right) ,x^{i}-\phi
^{1}\right\rangle -\frac{1}{2}\left\langle L\left( y^{i}-\phi ^{2}\right)
,y^{i}-\phi ^{2}\right\rangle -\frac{1}{2}\left\langle
Ru^{i},u^{i}\right\rangle .  \label{Hamiltonian function}
\end{eqnarray}%
Since $U$ is a closed convex set, then maximum principle reads as the
following local form
\begin{equation}
\left\langle \frac{\partial \mathcal{H}^{i}\left( t,p^{i,\ast },q^{i,\ast
},k^{i,\ast },x^{i,\ast },y^{i,\ast },z^{i,\ast },u^{i,\ast }\right) }{%
\partial u^{i}},u-u^{i,\ast }\right\rangle \leq 0,\text{ for }u\in U,\text{ }%
t\in \left[ 0,T\right] ,\text{ }P\text{-a.s.}
\label{convex maximum principle}
\end{equation}%
Hereafter, time argument is suppressed in case when no confusion occurs.
Noticing (\ref{Hamiltonian function}), then (\ref{convex maximum principle})
yields that
\begin{equation*}
\left\langle B^{T}q^{i,\ast }+K^{T}p^{i,\ast }+D^{T}k^{i,\ast }-Ru^{i,\ast
},u-u^{i,\ast }\right\rangle ,\text{ for all }u\in U,\text{ a.e. }t\in \left[
0,T\right] ,\text{ }P\text{-a.s.}
\end{equation*}%
or equivalently (noticing $R>0$),
\begin{equation}
\left\langle R^{\frac{1}{2}}\left[ R^{-1}\left( B^{T}q^{i,\ast
}+K^{T}p^{i,\ast }+D^{T}k^{i,\ast }\right) -u^{i,\ast }\right] ,R^{\frac{1}{2%
}}\left( u-u^{i,\ast }\right) \right\rangle ,
\label{optimal control condition}
\end{equation}%
for all $u\in U,$ a.e. $t\in \left[ 0,T\right] ,$ $P$-a.s$.$ As $R>0$, we
take the following norm on $U\subset \mathbb{R}^{m}$ (which is equivalent to
its Euclidean norm) $\Vert x\Vert _{R}^{2}=\left\langle \left\langle
x,x\right\rangle \right\rangle :=\left\langle R^{\frac{1}{2}}x,R^{\frac{1}{2}%
}x\right\rangle ,$ and by the well-known results of convex analysis, we
obtain that (\ref{optimal control condition}) is equivalent to
\begin{equation*}
u_{t}^{i,\ast }=\mathbf{P}_{U}\left[ R_{t}^{-1}\left( B_{t}^{T}q_{t}^{i,\ast
}+K_{t}^{T}p_{t}^{i,\ast }+D_{t}^{T}k_{t}^{i,\ast }\right) \right] ,
\end{equation*}%
where $\mathbf{P}_{U}(\cdot )$ is the projection mapping from $\mathbb{R}%
^{m} $ to its closed convex subset $\Gamma $ under the norm $\Vert \cdot
\Vert _{R}$. For more details, see Appendix. Hereafter, denote
\begin{equation*}
\varphi \left( t,p,q,k\right) =\mathbf{P}_{U}\left[ R_{t}^{-1}\left(
B_{t}^{T}q+K_{t}^{T}p+D_{t}^{T}k\right) \right] .
\end{equation*}%
Here, for simplicity, the dependence of $\varphi $ on time variable $t$ is
suppressed. The related Hamiltonian system becomes%
\begin{equation*}
\left\{
\begin{array}{rcl}
\mathrm{d}x_{t}^{i,\ast } & = & \left[ A_{t}x_{t}^{i,\ast }+B_{t}\varphi
\left( p_{t}^{i,\ast },q_{t}^{i,\ast },k_{t}^{i,\ast }\right) +F_{t}\phi
_{t}^{1}+b_{t}\right] \mathrm{d}t \\
&  & +\left[ D_{t}\varphi \left( p_{t}^{i,\ast },q_{t}^{i,\ast
},k_{t}^{i,\ast }\right) +\sigma _{t}\right] \mathrm{d}W_{t}^{i}, \\
\mathrm{d}y_{t}^{i,\ast } & = & -\left[ M_{t}x_{t}^{i,\ast
}+U_{t}y_{t}^{i,\ast }+H_{t}\phi _{t}^{1}+V_{t}\phi _{t}^{2}+K_{t}\varphi
\left( p_{t}^{i,\ast },q_{t}^{i,\ast },k_{t}^{i,\ast }\right) +f_{t}\right]
\mathrm{d}t+z_{t}^{i,\ast }\mathrm{d}W_{t}^{i} \\
\mathrm{d}p_{t}^{i,\ast } & = & \left[ U_{t}p_{t}^{i,\ast }-L_{t}\left(
y_{t}^{i,\ast }-\phi _{t}^{2}\right) \right] \mathrm{d}t, \\
\mathrm{d}q_{t}^{i,\ast } & = & \left[ -M_{t}p_{t}^{i,\ast
}-A_{t}q_{t}^{i,\ast }+Q_{t}\left( x_{t}^{i,\ast }-\phi _{t}^{1}\right) %
\right] \mathrm{d}t+k_{t}^{i,\ast }\mathrm{d}W_{t}^{i}, \\
x_{0}^{i,\ast } & = & x\in \mathbb{R}^{n},\text{ }y_{T}^{i,\ast }=\Phi
x_{T}^{i,\ast },\text{ }p_{0}^{i,\ast }=0,\text{ }q_{T}^{i,\ast }=\Phi
^{T}p_{T}^{i,\ast }-G\left( x_{T}^{i,\ast }-\phi _{T}^{1}\right) .%
\end{array}%
\right.
\end{equation*}%
After above preparations, it follows that
\begin{eqnarray}
\phi _{\cdot }^{1} &=&\lim_{N\rightarrow +\infty }\frac{1}{N}%
\sum_{i=1}^{N}x_{\cdot }^{i,\ast }=\mathbb{E}x_{\cdot }^{i,\ast },
\label{consi1} \\
\phi _{\cdot }^{2} &=&\lim_{N\rightarrow +\infty }\frac{1}{N}%
\sum_{i=1}^{N}y_{\cdot }^{i,\ast }=\mathbb{E}y_{\cdot }^{i,\ast }.
\label{consi2}
\end{eqnarray}%
Here, the first equality of \eqref{consi1} and \eqref{consi2} is due to the
consistency condition: the frozen term $\phi ^{1}$ and $\phi ^{2}$ should
equal to the average limit of all realized states $\left( x^{i,\ast
},y^{i,\ast }\right) ;$ the second equality is due to the law of large
numbers. Thus, by replacing $\left( \phi ^{1},\phi ^{2}\right) $ by $\left(
\mathbb{E}x^{i,\ast },\mathbb{E}y^{i,\ast }\right) $ in above Hamiltonian
system, we get the following system%
\begin{equation*}
\left\{
\begin{array}{rcl}
\mathrm{d}x_{t}^{i,\ast } & = & \left[ A_{t}x_{t}^{i,\ast }+B_{t}\varphi
\left( p_{t}^{i,\ast },q_{t}^{i,\ast },k_{t}^{i,\ast }\right) +F_{t}\mathbb{E%
}x_{t}^{i,\ast }+b_{t}\right] \mathrm{d}t \\
&  & +\left[ D_{t}\varphi \left( p_{t}^{i,\ast },q_{t}^{i,\ast
},k_{t}^{i,\ast }\right) +\sigma _{t}\right] \mathrm{d}W_{t}^{i}, \\
\mathrm{d}y_{t}^{i,\ast } & = & -\left[ M_{t}x_{t}^{i,\ast
}+U_{t}y_{t}^{i,\ast }+H_{t}\mathbb{E}x_{t}^{i,\ast }+V_{t}\mathbb{E}%
y_{t}^{i,\ast }+K_{t}\varphi \left( p_{t}^{i,\ast },q_{t}^{i,\ast
},k_{t}^{i,\ast }\right) +f_{t}\right] \mathrm{d}t+z_{t}^{i,\ast }\mathrm{d}%
W_{t}^{i}, \\
\mathrm{d}p_{t}^{i,\ast } & = & \left[ U_{t}p_{t}^{i,\ast }-L_{t}\left(
y_{t}^{i,\ast }-\mathbb{E}y_{t}^{i,\ast }\right) \right] \mathrm{d}t, \\
\mathrm{d}q_{t}^{i,\ast } & = & \left[ -M_{t}p_{t}^{i,\ast
}-A_{t}q_{t}^{i,\ast }+Q_{t}\left( x_{t}^{i,\ast }-\mathbb{E}x_{t}^{i,\ast
}\right) \right] \mathrm{d}t+k_{t}^{i,\ast }\mathrm{d}W_{t}^{i}, \\
x_{0}^{i,\ast } & = & =x\in \mathbb{R}^{n},\text{ }y_{T}^{i,\ast }=\Phi
x_{T}^{i,\ast },\text{ }p_{0}^{i}=0, \\
q_{T}^{i,\ast } & = & \Phi ^{T}p_{T}^{i,\ast }-G\left( x_{T}^{i,\ast }-%
\mathbb{E}x_{T}^{i,\ast }\right) .%
\end{array}%
\right.
\end{equation*}%
Clearly, all agents are statistically identical, therefore we can suppress
subscript \textquotedblleft $i$\textquotedblright\ and the following
consistency condition system appears for generic agent:%
\begin{equation}
\left\{
\begin{array}{rcl}
\mathrm{d}x & = & \left[ Ax+B\varphi \left( p,q,k\right) +F\mathbb{E}x+b%
\right] \mathrm{d}t+\left[ D\varphi \left( p,q,k\right) +\sigma \right]
\mathrm{d}W_{t}, \\
\mathrm{d}y & = & -\left[ Mx+Uy+H\mathbb{E}x+V\mathbb{E}y+K\varphi \left(
p,q,k\right) +f\right] \mathrm{d}t+z\mathrm{d}W_{t}, \\
\mathrm{d}p & = & \left[ Up-L\left( y-\mathbb{E}y\right) \right] \mathrm{d}t,
\\
\mathrm{d}q & = & \left[ -Mp-Aq+Q\left( x-\mathbb{E}x\right) \right] \mathrm{%
d}t+kW_{t}, \\
x_{0} & = & x\in \mathbb{R}^{n},\text{ }y_{T}=\Phi x_{T},\text{ }p_{0}=0, \\
q_{T} & = & \Phi ^{T}p_{T}-G\left( x_{T}-\mathbb{E}x_{T}\right) .%
\end{array}%
\right.  \label{cc}
\end{equation}%
Here, $W$ stands for a generic Brownian motion on $(\Omega ,\mathcal{F},P),$
and denote $\mathbb{F}^{W}$ the natural filtration generated by it and
augmented by all null-sets. $L_{\mathbb{F}^{W}}^{2},L_{\mathbb{F}^{W}}^{2,%
\mathcal{E}_{0}}$ are defined in the similar way with $L_{\mathbb{F}}^{2},L_{%
\mathbb{F}}^{2,\mathcal{E}_{0}}$ before. The system \eqref{cc} is a
nonlinear mean-field forward-backward SDE (MF-FBSDE) with projection
operator. It characterizes the state-average limit $\phi ^{1}=\mathbb{E}{x}%
,\phi ^{2}=\mathbb{E}{y}$ and MFG strategies $\bar{u}_{i}=\varphi (p,q,k)$
for a generic agent in the combined manner, which is totally different from
\cite{HHL, HHN}. As you may concern, we need to prove the above consistency
condition system admits a unique solution. We have the following uniqueness
and existence result.

\begin{remark}
It is necessary to point out that there should put a term $y_{0}^{i,\ast
}-y_{0}^{(N),\ast }$ in $(\ref{cost1})$. But we claim that after taking
expectation, it will disappear. Indeed, according to (\ref{optimal control
condition}), one has $\lim_{N\rightarrow +\infty }y_{0}^{(N),\ast
}=\lim_{N\rightarrow +\infty }\frac{1}{N}\sum_{i=1}^{N}y_{0}^{i,\ast }=%
\mathbb{E}y_{0}^{i,\ast }=y_{0}^{i,\ast }$. The first equality is just the
definition of $y_{0}^{(N),\ast },$ the second one is because of the law of
large numbers, whilst the last one is due to the fact that $y_{0}^{i,\ast }$
is an $\mathbb{F}_{0}^{i}$-measurable random vector; and therefore is
deterministic. Apparently, in contrast to Huang et al. \cite{HHL, HHN}, our
framework involves the state $\left( x^{i},y^{i}\right) .$
\end{remark}

\begin{theorem}
\label{th1}Assume that \emph{(A1) and (A2)} are in force. There exists a
unique adapted solution $(x,y,z,p,q,k)\in L_{\mathbb{F}^{W}}^{2}(0,T;\mathbb{%
R}^{n})\times L_{\mathbb{F}^{W}}^{2}(0,T;\mathbb{R}^{n})\times L_{\mathbb{F}%
^{W}}^{2}(0,T;\mathbb{R}^{n})\times L_{\mathbb{F}^{W}}^{2,\mathcal{E}%
_{0}}(0,T;\mathbb{R}^{n})\times L_{\mathbb{F}^{W}}^{2,\mathcal{E}_{0}}(0,T;%
\mathbb{R}^{n})\times L_{\mathbb{F}^{W}}^{2}(0,T;\mathbb{R}^{n})$ to system %
\eqref{cc}.
\end{theorem}

For simplicity, we put the proof of Theorem \ref{th1} in the Appendix \ref%
{app1}.

\section{ Main result}

\label{sect3}

In above sections, we can characterize the decentralized strategies $\{\bar{u%
}_{t}^{i},1\leq i\leq N\}$ of Problem (\textbf{CC}) through the auxiliary (%
\textbf{LCC}) and consistency condition system. For sake of presentation, we
alter the notation of consistency condition system to be $(\alpha ^{i},\beta
^{i},\gamma ^{i},\theta ^{i},\kappa ^{i},\gamma ^{i})$:%
\begin{equation}
\left\{
\begin{array}{rcl}
\mathrm{d}\alpha ^{i} & = & \left[ A\alpha ^{i}+B\varphi \left( \chi
^{i},\beta ^{i},\gamma ^{i}\right) +F\mathbb{E}\alpha ^{i}+b\right] \mathrm{d%
}t+\left[ D\varphi \left( \chi ^{i},\beta ^{i},\gamma ^{i}\right) +\sigma %
\right] \mathrm{d}W_{t}^{i}, \\
\mathrm{d}\theta ^{i} & = & -\left[ M\alpha ^{i}+U\theta ^{i}+H\mathbb{E}%
\alpha ^{i}+V\mathbb{E}\theta ^{i}+K\varphi \left( \chi ^{i},\beta
^{i},\gamma ^{i}\right) +f\right] \mathrm{d}t+\kappa ^{i}\mathrm{d}W_{t}^{i},
\\
\mathrm{d}\chi ^{i} & = & \left[ U\chi ^{i}-L\left( \theta ^{i}-\mathbb{E}%
\theta ^{i}\right) \right] \mathrm{d}t, \\
\mathrm{d}\beta ^{i} & = & \left[ -M\chi ^{i}-A\beta ^{i}+Q\left( \alpha
^{i}-\mathbb{E}\alpha ^{i}\right) \right] \mathrm{d}t+\gamma ^{i}W_{t}^{i},
\\
\alpha _{0}^{i} & = & x\in \mathbb{R}^{n},\text{ }\theta _{T}^{i}=\Phi
\alpha _{T}^{i},\text{ }\chi _{0}^{i}=0, \\
\beta _{T}^{i} & = & \Phi _{T}^{T}\chi ^{i}-G\left( \alpha _{T}^{i}-\mathbb{E%
}\alpha _{T}^{i}\right) .%
\end{array}%
\right.  \label{lastadded}
\end{equation}%
Now, we are in position to verify the $\epsilon $-Nash equilibrium of them.
To this end, let us first present the definition of $\epsilon $-Nash
equilibrium.

\begin{definition}
\label{d1} A set of strategies, $\bar u_t^i\in \mathcal{U}^{c}_{ad}$, $1\leq
i\leq N$, for $N$ agents, is called to satisfy an $\epsilon$-Nash
equilibrium with respect to costs $\mathcal{J}^i,\ 1\leq i\leq N,$ if there
exists $\epsilon=\epsilon(N)\geq0, \displaystyle\lim_{N \rightarrow
+\infty}\epsilon(N)=0,$ such that for any $1\leq i\leq N$, we have
\begin{equation}  \label{epsilon}
\mathcal{J}^i(\bar{u}_t^i,\bar{u}_t^{-i})\leq \mathcal{J}^i(u_t^i,\bar{u}%
_t^{-i})+\epsilon,
\end{equation}
when any alternative strategy $u^i\in \mathcal{U}^{c}_{ad}$ is applied by $%
\mathcal{A}_i$.
\end{definition}

\begin{remark}
If $\epsilon=0$, then Definition \ref{d1} is reduced to the usual exact Nash
equilibrium.
\end{remark}

Now, we give the main result of this paper and its proof will be shown step
by step.

\begin{theorem}
\label{Nash equilibrium theorem} Assume that \emph{(A1) and (A2) }are in
force. Then, $(\bar{u}_{1},\bar{u}_{2},\cdots ,\bar{u}_{N})$ is an $\epsilon
$-Nash equilibrium of \emph{\textbf{Problem (CC)}}.
\end{theorem}

In order to prove the main Theorem \ref{Nash equilibrium theorem}, we needs
several lemmas which are presented later. For agent $\mathcal{A}_{i},$
recall that its decentralized open-loop optimal strategy is $\bar{u}%
_{i}=\varphi \left( \chi ^{i},\beta ^{i},\gamma ^{i}\right) $. The
decentralized state $\left( \breve{x}_{t}^{i},\breve{y}^{i},\breve{z}%
^{i}\right) ,$ is%
\begin{equation}
\left\{
\begin{array}{rcl}
\mathrm{d}\breve{x}^{i} & = & \left[ A\breve{x}^{i}+B\varphi \left( \chi
^{i},\beta ^{i},\gamma ^{i}\right) +F\breve{x}^{(N)}+b\right] \mathrm{d}t+%
\left[ D\varphi \left( \chi ^{i},\beta ^{i},\gamma ^{i}\right) +\sigma %
\right] \mathrm{d}W_{t}^{i}, \\
\mathrm{d}\breve{y}^{i} & = & -\left[ M\breve{x}^{i}+U\breve{y}^{i}+H\breve{x%
}^{(N)}+V\breve{y}^{(N)}+K\varphi \left( \chi ^{i},\beta ^{i},\gamma
^{i}\right) +f\right] \mathrm{d}t+\breve{z}^{i}\mathrm{d}W_{t}^{i}, \\
\mathrm{d}\alpha ^{i} & = & \left[ A\alpha ^{i}+B\varphi \left( \chi
^{i},\beta ^{i},\gamma ^{i}\right) +F\mathbb{E}\alpha ^{i}+b\right] \mathrm{d%
}t+\left[ D\varphi \left( \chi ^{i},\beta ^{i},\gamma ^{i}\right) +\sigma %
\right] \mathrm{d}W_{t}^{i}, \\
\mathrm{d}\theta ^{i} & = & -\left[ M\alpha ^{i}+U\theta ^{i}+H\mathbb{E}%
\alpha ^{i}+V\mathbb{E}\theta ^{i}+K\varphi \left( \chi ^{i},\beta
^{i},\gamma ^{i}\right) +f\right] \mathrm{d}t+\kappa ^{i}\mathrm{d}W_{t}^{i},
\\
\mathrm{d}\chi ^{i} & = & \left[ U\chi ^{i}-L\left( \theta ^{i}-\mathbb{E}%
\theta ^{i}\right) \right] \mathrm{d}t, \\
\mathrm{d}\beta ^{i} & = & \left[ -M\chi ^{i}-A\beta ^{i}+Q\left( \alpha
^{i}-\mathbb{E}\alpha ^{i}\right) \right] \mathrm{d}t+\gamma ^{i}W_{t}^{i},
\\
\breve{x}_{0}^{i} & = & \alpha _{0}^{i}=x,\text{ }\breve{y}_{T}^{i}=\Phi
\breve{x}_{T}^{i}, \\
\theta _{T}^{i} & = & \Phi \alpha _{T}^{i},\text{ }\chi _{0}^{i}=-\Psi
\left( \theta _{0}^{i}-\mathbb{E}\theta _{0}^{i}\right) , \\
\beta _{T}^{i} & = & \Phi _{T}^{T}\chi ^{i}-G\left( \alpha _{T}^{i}-\mathbb{E%
}\alpha _{T}^{i}\right) ,%
\end{array}%
\right.  \label{decentralized state}
\end{equation}%
where $\breve{x}^{(N)}=\frac{1}{N}\sum_{i=1}^{N}\breve{x}^{i}$ and $\breve{y}%
^{(N)}=\frac{1}{N}\sum_{i=1}^{N}\breve{y}^{i}$ Note that $(\alpha ^{i},\beta
^{i},\gamma ^{i},\theta ^{i},\kappa ^{i},\gamma ^{i})$ satisfies (\ref%
{lastadded}).

For each $1\leq i\leq N$, the monotonic fully coupled FBSDEs (\ref{lastadded}%
) has a unique solution $(\alpha ^{i},\beta ^{i},\gamma ^{i})\in L_{\mathbb{F%
}^{i}}^{2}(0,T;\mathbb{R}^{n})\times L_{\mathbb{F}^{i}}^{2}(0,T;\mathbb{R}%
^{n})\times L_{\mathbb{F}^{i}}^{2}(0,T;\mathbb{R}^{n})$. Thus, the system of
all first equation of (\ref{decentralized state}), $1\leq i\leq N$, has also
a unique solution $\left( (\breve{x}^{i})_{i},(\breve{y}^{i})_{i},(\breve{%
\kappa}^{i})_{i}\right) \in (L_{\mathbb{F}}^{2}(0,T;\mathbb{R}%
^{n}))^{\otimes N}\times (L_{\mathbb{F}}^{2}(0,T;\mathbb{R}^{n}))^{\otimes
N}\times (L_{\mathbb{F}}^{2}(0,T;\mathbb{R}^{n}))^{\otimes N}$, where $%
\otimes N$ denotes the $n$-tuple Cartesian product. Moreover, since $%
\{W_{i}\}_{i=1}^{N}$ is $N$-dimensional Brownian motions whose components
are independent and identically distributed, we have $(\alpha ^{i},\beta
^{i},\gamma ^{i},\theta ^{i},\kappa ^{i},\gamma ^{i}),1\leq i\leq N$ are
independent and identically distributed.

\begin{lemma}
\label{estimatesdebsde}If \emph{(A1) and (A2)} hold, then
\begin{eqnarray*}
\mathbb{E}\sup_{0\leq t\leq T}\Big|\breve{x}^{(N)}(t)\Big|^{2} &\leq &C_{2},
\\
\mathbb{E}\sup_{0\leq t\leq T}\Big|\breve{y}^{(N)}(t)\Big|^{2} &\leq &C_{2}.
\end{eqnarray*}
\end{lemma}

The proof of Lemma \ref{estimatesdebsde} is classical by virtue of B-D-G
inequality and Schwarz inequality, so we omit it.

\begin{lemma}
\label{lemma for xN and m}If \emph{(A1) and (A2)} hold, then
\begin{eqnarray}
\mathbb{E}\sup_{0\leq t\leq T}\Big|\breve{x}^{(N)}(t)-\mathbb{E}\alpha
^{i}(t)\Big|^{2} &=&O\Big(\frac{1}{N}\Big),  \label{sde1} \\
\mathbb{E}\sup_{0\leq t\leq T}\Big|\breve{y}^{(N)}(t)-\mathbb{E}\theta
^{i}(t)\Big|^{2} &=&O\Big(\frac{1}{N}\Big).  \label{bsde1}
\end{eqnarray}
\end{lemma}

\paragraph{Proof.}

Let us add up both sides of the first and second equation of (\ref%
{decentralized state}) with respect to all $1\leq i\leq N$ and multiply $%
\frac{1}{N}$, we obtain (recall that $\breve{x}^{(N)}=\frac{1}{N}%
\sum_{i=1}^{N}\breve{x}^{i}$, $\breve{y}^{(N)}=\frac{1}{N}\sum_{i=1}^{N}%
\breve{y}^{i}$ and $\breve{z}^{(N)}=\frac{1}{N}\sum_{i=1}^{N}\breve{z}^{i}$
\begin{equation}
\left\{
\begin{array}{rl}
\mathrm{d}\breve{x}^{(N)}= & \left[ A\breve{x}^{(N)}+\frac{1}{N}%
\sum_{i=1}^{N}B\varphi \left( \chi ^{i},\beta ^{i},\gamma ^{i}\right) +F%
\breve{x}^{(N)}+b\right] \mathrm{d}t \\
& +\displaystyle\frac{1}{N}\sum_{i=1}^{N}\left[ D\varphi \left( \chi
^{i},\beta ^{i},\gamma ^{i}\right) +\sigma \right] \mathrm{d}W_{t}^{i}, \\
\mathrm{d}\breve{y}^{(N)}= & -\left[ M\breve{x}^{(N)}+U\breve{y}^{(N)}+H%
\breve{x}^{(N)}+V\breve{y}^{(N)}+\frac{1}{N}\sum_{i=1}^{N}K\varphi \left(
\chi ^{i},\beta ^{i},\gamma ^{i}\right) +f\right] \mathrm{d}t \\
& +\displaystyle\frac{1}{N}\sum_{i=1}^{N}\breve{z}^{i}\mathrm{d}W_{t}^{i},
\\
\breve{x}_{0}^{(N)}= & x,\quad \breve{y}_{T}^{(N)}=\Phi \breve{x}_{T}^{(N)}.%
\end{array}%
\right.  \label{ee1}
\end{equation}%
On the other hand, by taking the expectation on both sides of the second
equation of (\ref{decentralized state}), it follows from Fubini's theorem
that $\mathbb{E}\alpha ^{i}$ satisfies the following equation:
\begin{equation}
\left\{
\begin{array}{l}
\mathrm{d}\left( \mathbb{E}\alpha ^{i}\right) =\left[ A\mathbb{E}\alpha
^{i}+B\mathbb{E}\varphi \left( \chi ^{i},\beta ^{i},\gamma ^{i}\right) +F%
\mathbb{E}\alpha ^{i}+b\right] \mathrm{d}t, \\
\mathrm{d}\left( \mathbb{E}\theta ^{i}\right) =-\left[ M\mathbb{E}\alpha
^{i}+U\mathbb{E}\theta ^{i}+H\mathbb{E}\alpha ^{i}+V\mathbb{E}\theta ^{i}+K%
\mathbb{E}\varphi \left( \chi ^{i},\beta ^{i},\gamma ^{i}\right) +f\right]
\mathrm{d}t, \\
\mathbb{E}\alpha _{0}^{i}=x,\quad \mathbb{E}\theta _{T}^{i}=\Phi \mathbb{E}%
\alpha _{T}^{i}.%
\end{array}%
\right.  \label{ee2}
\end{equation}%
From (\ref{ee1}) and (\ref{ee2}), by denoting
\begin{eqnarray}
\Delta _{t}^{1} &:&=\breve{x}^{(N)}(t)-\mathbb{E}\alpha ^{i}(t),
\label{insde1} \\
\Delta _{t}^{2} &:&=\breve{y}^{(N)}(t)-\mathbb{E}\theta ^{i}(t),
\label{inbsde3}
\end{eqnarray}%
we have
\begin{equation*}
\left\{
\begin{array}{l}
\mathrm{d}\Delta ^{1}=\left[ A\Delta ^{1}+\frac{1}{N}\sum_{i=1}^{N}B\varphi
\left( \chi ^{i},\beta ^{i},\gamma ^{i}\right) -B\mathbb{E}\varphi \left(
\chi ^{i},\beta ^{i},\gamma ^{i}\right) +F\Delta ^{1}\right] \mathrm{d}t \\
\qquad \qquad +\displaystyle\frac{1}{N}\sum_{i=1}^{N}\left[ D\varphi \left(
\chi ^{i},\beta ^{i},\gamma ^{i}\right) +\sigma \right] \mathrm{d}W_{t}^{i},
\\
\mathrm{d}\Delta ^{2}=-\Big [M\Delta ^{1}+U\Delta ^{2}+H\Delta ^{1}+V\Delta
^{2} \\
\qquad \qquad +\displaystyle\frac{1}{N}\sum_{i=1}^{N}K\varphi \left( \chi
^{i},\beta ^{i},\gamma ^{i}\right) -K\mathbb{E}\varphi \left( \chi
^{i},\beta ^{i},\gamma ^{i}\right) \Big ]\mathrm{d}t+\frac{1}{N}%
\sum_{i=1}^{N}\breve{z}^{i}\mathrm{d}W_{t}^{i}, \\
\Delta _{0}^{1}=0,\quad \Delta _{T}^{2}=\Phi \Delta _{T}^{1},%
\end{array}%
\right.
\end{equation*}%
and the inequality $(x+y)^{2}\leq 2x^{2}+2y^{2}$ yields that, for any $t\in
\lbrack 0,T]$,%
\begin{eqnarray*}
\mathbb{E}\left[ \sup_{0\leq s\leq t}\left\vert \Delta _{s}^{1}\right\vert
^{2}\right] &\leq &2\mathbb{E}\sup_{0\leq s\leq t}\Bigg |\int_{0}^{s}\Big [%
\left( A+F\right) \Delta _{s}^{1} \\
&&+\frac{1}{N}\sum_{i=1}^{N}B\varphi \left( \chi ^{i},\beta ^{i},\gamma
^{i}\right) -B\mathbb{E}\varphi \left( \chi ^{i},\beta ^{i},\gamma
^{i}\right) \Big ]\mathrm{d}r\Bigg |^{2} \\
&&+2\mathbb{E}\sup_{0\leq s\leq t}\Bigg |\frac{1}{N}\sum_{i=1}^{N}%
\int_{0}^{s}\left[ D\varphi \left( \chi ^{i},\beta ^{i},\gamma ^{i}\right)
+\sigma \right] \mathrm{d}W_{t}^{i}\Bigg |^{2}.
\end{eqnarray*}%
From the well-known Cauchy-Schwartz inequality and the B-D-G inequality, we
obtain that there exists a constant $C_{0}$ independent of $N$ (which may
vary line by line) such that%
\begin{eqnarray}
\mathbb{E}\left[ \sup_{0\leq s\leq t}\left\vert \Delta _{s}^{1}\right\vert
^{2}\right] &\leq &C_{0}\mathbb{E}\sup_{0\leq s\leq t}\Bigg |\int_{0}^{s}%
\Bigg [\left\vert \Delta _{s}^{1}\right\vert ^{2}+\left\vert \frac{1}{N}%
\sum_{i=1}^{N}B\varphi \left( \chi ^{i},\beta ^{i},\gamma ^{i}\right) -B%
\mathbb{E}\varphi \left( \chi ^{i},\beta ^{i},\gamma ^{i}\right) \right\vert
^{2}\Bigg ]\mathrm{d}r\Bigg |^{2}  \notag \\
&&+\frac{C_{0}}{N}\mathbb{E}\Bigg |\sum_{i=1}^{N}\int_{0}^{s}\left\vert
D\varphi \left( \chi ^{i},\beta ^{i},\gamma ^{i}\right) +\sigma \right\vert
^{2}\mathrm{d}t\Bigg |^{2}.  \label{ee3}
\end{eqnarray}%
Since $(\chi ^{i},\beta ^{i},\gamma ^{i}),1\leq i\leq N$ are independent
identically distributed, for each fixed $s\in \lbrack 0,T]$, let us denote
that $\mu (s)=\mathbb{E}\varphi (\chi ^{i},\beta ^{i},\gamma ^{i}))$ (note
that $\mu $ does not depend on $i$), we have%
\begin{eqnarray*}
&&\mathbb{E}\left\vert \frac{1}{N}\sum_{i=1}^{N}\varphi \left( \chi
^{i},\beta ^{i},\gamma ^{i}\right) -\mathbb{E}\varphi \left( \chi ^{i},\beta
^{i},\gamma ^{i}\right) \right\vert ^{2} \\
&=&\frac{1}{N^{2}}\mathbb{E}\left\vert \sum_{i=1}^{N}\varphi \left( \chi
^{i},\beta ^{i},\gamma ^{i}\right) -\mu _{s}\right\vert ^{2} \\
&=&\frac{1}{N^{2}}\mathbb{E}\sum_{i=1}^{N}\left\vert \varphi \left( \chi
^{i},\beta ^{i},\gamma ^{i}\right) -\mu _{s}\right\vert ^{2} \\
&&+\frac{2}{N^{2}}\mathbb{E}\left[ \sum_{i=1,j=1,j\neq i}^{N}\left\langle
\varphi \left( \chi ^{i},\beta ^{i},\gamma ^{i}\right) -\mu _{s},\varphi
\left( \chi ^{j},\beta ^{j},\gamma ^{j}\right) -\mu _{s}\right\rangle \right]
.
\end{eqnarray*}%
Since $(\chi ^{i},\beta ^{i},\gamma ^{i}),1\leq i\leq N$ are independent, we
have%
\begin{eqnarray*}
&&\frac{2}{N^{2}}\mathbb{E}\left[ \sum_{i=1,j=1,j\neq i}^{N}\left\langle
\varphi \left( \chi ^{i},\beta ^{i},\gamma ^{i}\right) -\mu _{s},\varphi
\left( \chi ^{j},\beta ^{j},\gamma ^{j}\right) -\mu _{s}\right\rangle \right]
\\
&=&\frac{2}{N^{2}}\sum_{i=1,j=1,j\neq i}^{N}\left\langle \mathbb{E}\varphi
\left( \chi ^{i},\beta ^{i},\gamma ^{i}\right) -\mu _{s},\mathbb{E}\varphi
\left( \chi ^{j},\beta ^{j},\gamma ^{j}\right) -\mu _{s}\right\rangle =0.
\end{eqnarray*}%
Then, due to the fact that $(\chi ^{i},\beta ^{i},\gamma ^{i}),1\leq i\leq N$
are identically distributed, there exists a constant $C_{0}$ independent of $%
N$ such that%
\begin{eqnarray*}
&&\mathbb{E}\int_{0}^{t}\left\vert \frac{1}{N}\sum_{i=1}^{N}B\varphi \left(
\chi ^{i},\beta ^{i},\gamma ^{i}\right) -B\mathbb{E}\varphi \left( \chi
^{i},\beta ^{i},\gamma ^{i}\right) \right\vert ^{2}\mathrm{d}r \\
&\leq &C_{0}\int_{0}^{t}\mathbb{E}\left\vert \frac{1}{N}\sum_{i=1}^{N}%
\varphi \left( \chi ^{i},\beta ^{i},\gamma ^{i}\right) -\mu _{s}\right\vert
^{2}\mathrm{d}r \\
&\leq &\frac{C_{0}}{N^{2}}\int_{0}^{t}\mathbb{E}\sum_{i=1}^{N}\left\vert
\varphi \left( \chi ^{i},\beta ^{i},\gamma ^{i}\right) -\mu _{s}\right\vert
^{2}\mathrm{d}r \\
&=&\frac{C_{0}}{N}\int_{0}^{t}\mathbb{E}\left\vert \varphi \left( \chi
^{i},\beta ^{i},\gamma ^{i}\right) -\mu _{s}\right\vert ^{2}\mathrm{d}r \\
&=&O\left( \frac{1}{N}\right) ,
\end{eqnarray*}%
where the last equality comes from the fact that $\varphi (\chi ^{i},\beta
^{i},\gamma ^{i})\in L_{\mathcal{F}^{i}}^{2}(0,T;\Gamma )$.

We proceed the second term of (\ref{ee3}), using the fact that $(\chi
^{i},\beta ^{i},\gamma ^{i})$ are identically distributed as follows:
\begin{equation*}
\frac{C_{0}}{N^{2}}\mathbb{E}\Bigg[\sum_{i=1}^{N}\int_{0}^{t}\left\vert
\!D\varphi (\chi ^{i},\beta ^{i},\gamma ^{i})+\sigma (s)\right\vert ^{2}ds%
\Bigg]=O\Big(\frac{1}{N}\Big).
\end{equation*}%
Moreover, we obtain from (\ref{ee3}) that
\begin{equation*}
\mathbb{E}\sup_{0\leq s\leq t}|\Delta _{s}^{1}|^{2}\leq C_{0}\mathbb{E}%
\int_{0}^{t}|\Delta _{s}^{1}|^{2}+O\Big(\frac{1}{N}\Big),\text{ for any }%
t\in \lbrack 0,T].
\end{equation*}%
Consequently, by virtue of Gronwall's inequality, we get the first estimate (%
\ref{insde1}).

We now handle the estimates (\ref{inbsde3}). Applying It\^{o}'s formula
again, we have%
\begin{eqnarray*}
\left\vert \Delta _{t}^{2}\right\vert ^{2}+\frac{1}{N}\sum_{i=1}^{N}%
\int_{t}^{T}\left\vert \breve{z}_{s}^{i}\right\vert ^{2}\mathrm{d}s
&=&\left\vert \Delta _{T}^{2}\right\vert ^{2}+2\int_{t}^{T}\Big <\Delta
^{2},M\Delta ^{1}+U\Delta ^{2}+H\Delta ^{1}+V\Delta ^{2} \\
&&+\frac{1}{N}\sum_{i=1}^{N}K\varphi \left( \chi ^{i},\beta ^{i},\gamma
^{i}\right) -K\mathbb{E}\varphi \left( \chi ^{i},\beta ^{i},\gamma
^{i}\right) \Big >\mathrm{d}s \\
&&-2\frac{1}{N}\sum_{i=1}^{N}\int_{t}^{T}\Big <\Delta ^{2},\breve{z}^{i}\Big
>\mathrm{d}W_{s}^{i}.
\end{eqnarray*}%
Using B-D-G inequalities, we show that there exists a constant $C_{1},$
modifying $C_{1}$ if necessary,
\begin{eqnarray*}
&&\mathbb{E}\left[ \sup_{0\leq t\leq T}\left\vert \Delta _{t}^{2}\right\vert
^{2}+\frac{1}{N}\sum_{i=1}^{N}\int_{t}^{T}\left\vert \breve{z}%
^{i}\right\vert ^{2}\mathrm{d}s\right] \\
&\leq &\mathbb{E}\left\vert \Phi \Delta _{T}^{1}\right\vert ^{2}+C_{1}%
\mathbb{E}\Bigg [\int_{0}^{T}\Big <\Delta ^{2},M\Delta ^{1}\Big >+\Big <%
\Delta ^{2},U\Delta ^{2}\Big > \\
&&+\Big <\Delta ^{2},H\Delta ^{1}\Big >+\Big <\Delta ^{2},V\Delta ^{2}\Big >
\\
&&+\Big <\Delta ^{2},\frac{1}{N}\sum_{i=1}^{N}K\varphi \left( \chi
^{i},\beta ^{i},\gamma ^{i}\right) -K\mathbb{E}\varphi \left( \chi
^{i},\beta ^{i},\gamma ^{i}\right) \Big >\mathrm{d}s \\
&&+\frac{1}{N}\sum_{i=1}^{N}\left( \int_{0}^{T}\left\vert \Delta
^{2}\right\vert ^{2}\left\vert \breve{z}^{i}\right\vert ^{2}\mathrm{d}%
s\right) ^{\frac{1}{2}}\Bigg ] \\
&\leq &\mathbb{E}\left\vert \Phi \Delta _{T}^{1}\right\vert ^{2}+C_{1}%
\mathbb{E}\Bigg [\int_{0}^{T}\left\vert \Delta ^{2}\right\vert \Big (%
\left\vert M\Delta ^{1}\right\vert +\left\vert U\Delta ^{2}\right\vert \\
&&+\left\vert H\Delta ^{1}\right\vert +\left\vert V\Delta ^{2}\right\vert
+\left\vert \frac{1}{N}\sum_{i=1}^{N}K\varphi \left( \chi ^{i},\beta
^{i},\gamma ^{i}\right) -K\mathbb{E}\varphi \left( \chi ^{i},\beta
^{i},\gamma ^{i}\right) \right\vert \Big )\mathrm{d}s \\
&&+\frac{1}{N}\sum_{i=1}^{N}\left( \int_{0}^{T}\left\vert \Delta
^{2}\right\vert ^{2}\left\vert \breve{z}^{i}\right\vert ^{2}\mathrm{d}%
s\right) ^{\frac{1}{2}}\Bigg ].
\end{eqnarray*}%
Employing the classical Cauchy-Schwarz inequality and Gronwall's inequality
with estimation (\ref{sde1}), we get (\ref{bsde1}). \hfill $\Box $

\begin{lemma}
\label{lemma for x and xi}Assume that \emph{(A1) and (A2)} are in force.
Then, we have
\begin{eqnarray}
\sup_{1\leq i\leq N}\mathbb{E}\left[ \sup_{0\leq t\leq T}\Big|\breve{x}%
_{t}^{i}-\alpha _{t}^{i}\Big|^{2}\right] &=&O\Big(\frac{1}{N}\Big),
\label{sde2} \\
\sup_{1\leq i\leq N}\mathbb{E}\left[ \sup_{0\leq t\leq T}\Big|\breve{y}%
_{t}^{i}-\theta _{t}^{i}\Big|^{2}\right] &=&O\Big(\frac{1}{N}\Big).
\label{bsde2}
\end{eqnarray}
\end{lemma}

\paragraph{Proof.}

From (\ref{decentralized state}) and (\ref{lastadded}), we have that%
\begin{equation}
\left\{
\begin{array}{l}
\mathrm{d}\breve{x}^{i}=\left[ A\breve{x}^{i}+B\varphi \left( \chi
^{i},\beta ^{i},\gamma ^{i}\right) +F\breve{x}^{(N)}+b\right] \mathrm{d}t+%
\left[ D\varphi \left( \chi ^{i},\beta ^{i},\gamma ^{i}\right) +\sigma %
\right] \mathrm{d}W_{t}^{i}, \\
\mathrm{d}\breve{y}^{i}=-\left[ M\breve{x}^{i}+U\breve{y}^{i}+H\breve{x}%
^{(N)}+V\breve{y}^{(N)}+K\varphi \left( \chi ^{i},\beta ^{i},\gamma
^{i}\right) +f\right] \mathrm{d}t+\breve{z}^{i}\mathrm{d}W_{t}^{i}, \\
\mathrm{d}\alpha ^{i}=\left[ A\alpha ^{i}+B\varphi \left( \chi ^{i},\beta
^{i},\gamma ^{i}\right) +F\mathbb{E}\alpha ^{i}+b\right] \mathrm{d}t+\left[
D\varphi \left( \chi ^{i},\beta ^{i},\gamma ^{i}\right) +\sigma \right]
\mathrm{d}W_{t}^{i}, \\
\mathrm{d}\theta ^{i}=-\left[ M\alpha ^{i}+U\theta ^{i}+H\mathbb{E}\alpha
^{i}+V\mathbb{E}\theta ^{i}+K\varphi \left( \chi ^{i},\beta ^{i},\gamma
^{i}\right) +f\right] \mathrm{d}t+\kappa ^{i}\mathrm{d}W_{t}^{i}, \\
\breve{x}_{0}^{i}=\alpha _{0}^{i}=x,\text{ }\breve{y}_{T}^{i}=\Phi \breve{x}%
_{T}^{i},\text{ }\theta _{T}^{i}=\Phi \alpha _{T}^{i},\text{ }%
\end{array}%
\right.  \label{ee4}
\end{equation}%
where $(\beta ^{i},\chi ^{i},\gamma ^{i})$ is the unique solution to the
following FBSDEs:%
\begin{equation*}
\left\{
\begin{array}{l}
\mathrm{d}\chi ^{i}=\left[ U\chi ^{i}-L\left( \theta ^{i}-\mathbb{E}\theta
^{i}\right) \right] \mathrm{d}t, \\
\mathrm{d}\beta ^{i}=\left[ -M\chi ^{i}-A\beta ^{i}+Q\left( \alpha ^{i}-%
\mathbb{E}\alpha ^{i}\right) \right] \mathrm{d}t+\gamma ^{i}W_{t}^{i}, \\
\breve{x}_{0}^{i}=x,\text{ }\chi _{0}^{i}=0, \\
\beta _{T}^{i}=\Phi _{T}^{T}\chi ^{i}-G\left( \alpha _{T}^{i}-\mathbb{E}%
\alpha _{T}^{i}\right) .%
\end{array}%
\right.
\end{equation*}%
From (\ref{ee4}), we have
\begin{equation*}
\left\{ \begin{aligned}
d(\breve{x}^{i}-\alpha^{i})&=\Big[A(\breve{x}^{i}\!-\!\alpha^{i})
\!+\!F(\breve{x}^{(N)}\!-\!\mathbb{E}\alpha^{i})\Big]dt, \\
\breve{x}^{i}(0)-\bar{x}^{i}(0)&=0. \end{aligned}\right.
\end{equation*}%
The classical estimate for the SDE yields that
\begin{equation*}
\mathbb{E}\sup_{0\leq t\leq T}\Big|\breve{x}_{t}^{i}-\alpha _{t}^{i}\Big|%
^{2}\leq C_{0}\mathbb{E}\int_{0}^{T}\left\vert \breve{x}_{s}^{(N)}-\mathbb{E}%
\alpha _{s}^{i}\right\vert ^{2}\mathrm{d}s,
\end{equation*}%
where $C_{0}$ is a constant independent of $N$. Noticing (\ref{sde1}) of
Lemma \ref{lemma for xN and m}, we obtain (\ref{sde2}). We consider%
\begin{equation*}
\left\{
\begin{array}{l}
\mathrm{d}\left( \breve{y}_{t}^{i}-\theta _{t}^{i}\right) =-\Big [M\left(
\breve{x}^{i}-\alpha ^{i}\right) +U\left( \breve{y}_{t}^{i}-\theta
_{t}^{i}\right) \\
\qquad \qquad +F\left( \breve{x}^{(N)}-\mathbb{E}\alpha ^{i}\right) +V\left(
\breve{y}^{(N)}-\mathbb{E}\theta ^{i}\right) \Big ]\mathrm{d}t+\left( \breve{%
z}^{i}-\kappa ^{i}\right) \mathrm{d}W_{t}^{i}, \\
\breve{y}_{T}^{i}-\theta _{T}^{i}=\Phi \left( \breve{x}_{T}^{i}-\alpha
_{T}^{i}\right) .%
\end{array}%
\right.
\end{equation*}%
By classical estimation for BSDE, we have
\begin{eqnarray*}
&&\sup_{1\leq i\leq N}\mathbb{E}\left[ \sup_{0\leq t\leq T}\Big|\breve{y}%
_{t}^{i}-\theta _{t}^{i}\Big|^{2}\right] +\mathbb{E}\int_{0}^{T}\left\vert
\breve{z}^{i}-\kappa ^{i}\right\vert ^{2}\mathrm{d}s \\
&\leq &C_{0}\mathbb{E}\int_{0}^{T}\Bigg (\Big|\breve{x}^{(N)}-\mathbb{E}%
\alpha ^{i}\Big|^{2}+\Big|\breve{y}^{(N)}-\mathbb{E}\theta ^{i}\Big|^{2} \\
&&+\Big|\breve{x}^{i}-\alpha ^{i}\Big|^{2}+\Big|\breve{y}_{t}^{i}-\theta
_{t}^{i}\Big|^{2}\Bigg )\mathrm{d}s.
\end{eqnarray*}%
where $C_{0}$ is a constant independent of $N$. By Gronwall's inequality, we
get the desired result.\hfill $\Box $

\begin{lemma}
\label{first lemma for cost} For all $1\leq i\leq N$, we have
\begin{equation*}
\Big|\mathcal{J}_{i}(\bar{u}^{i},\bar{u}_{-i})-J_{i}(\bar{u}^{i})\Big|=O\Big(%
\frac{1}{\sqrt{N}}\Big).
\end{equation*}
\end{lemma}

\paragraph{Proof.}

From the definition of (\ref{cost1}), (\ref{cost2}) and (\ref{decentralized
state}), we have%
\begin{eqnarray*}
\mathcal{J}_{i}\left( \bar{u}^{i},\bar{u}_{-i}\right) &=&\frac{1}{2}\mathbb{E%
}\Bigg [\int_{0}^{T}\Big [\left\langle Q_{t}\left( \breve{x}_{t}^{i}-\breve{x%
}_{t}^{(N)}\right) ,\breve{x}_{t}^{i}-\breve{x}_{t}^{(N)}\right\rangle
+\left\langle L_{t}\left( \breve{y}_{t}^{i}-\breve{y}_{t}^{(N)}\right) ,%
\breve{y}_{t}^{i}-\breve{y}_{t}^{(N)}\right\rangle \\
&&+\left\langle R_{t}\bar{u}_{t}^{i},\bar{u}_{t}^{i}\right\rangle \Big ]%
\mathrm{d}t+\left\langle G\left( \breve{x}_{T}^{i}-\breve{x}%
_{T}^{(N)}\right) ,\breve{x}_{T}^{i}-\breve{x}_{T}^{(N)}\right\rangle \Bigg ]
\end{eqnarray*}%
and
\begin{eqnarray*}
J_{i}\left( \bar{u}^{i}\right) &=&\frac{1}{2}\mathbb{E}\Bigg [\int_{0}^{T}%
\Big [\left\langle Q_{t}\left( \alpha _{t}^{i}-\mathbb{E}\alpha
_{t}^{i}\right) ,\alpha _{t}^{i}-\mathbb{E}\alpha _{t}^{i}\right\rangle
+\left\langle L_{t}\left( \theta _{t}^{i}-\mathbb{E}\theta _{t}^{i}\right)
,\theta _{t}^{i}-\mathbb{E}\theta _{t}^{i}\right\rangle \\
&&+\left\langle R_{t}\bar{u}_{t}^{i},\bar{u}_{t}^{i}\right\rangle \Big ]%
\mathrm{d}t+\left\langle G\left( \alpha _{T}^{i}-\mathbb{E}\alpha
_{T}^{i}\right) ,\alpha _{T}^{i}-\mathbb{E}\alpha _{T}^{i}\right\rangle %
\Bigg ],
\end{eqnarray*}%
then
\begin{eqnarray}
\mathcal{J}_{i}\left( \bar{u}^{i},\bar{u}_{-i}\right) -J_{i}\left( \bar{u}%
^{i}\right) &=&\frac{1}{2}\mathbb{E}\Bigg [\int_{0}^{T}\Big [\left\langle
Q_{t}\left( \breve{x}_{t}^{i}-\breve{x}_{t}^{(N)}\right) ,\breve{x}_{t}^{i}-%
\breve{x}_{t}^{(N)}\right\rangle -\left\langle Q_{t}\left( \alpha _{t}^{i}-%
\mathbb{E}\alpha _{t}^{i}\right) ,\alpha _{t}^{i}-\mathbb{E}\alpha
_{t}^{i}\right\rangle  \notag \\
&&+\left\langle L_{t}\left( \breve{y}_{t}^{i}-\breve{y}_{t}^{(N)}\right) ,%
\breve{y}_{t}^{i}-\breve{y}_{t}^{(N)}\right\rangle -\left\langle L_{t}\left(
\theta _{t}^{i}-\mathbb{E}\theta _{t}^{i}\right) ,\theta _{t}^{i}-\mathbb{E}%
\theta _{t}^{i}\right\rangle \mathrm{d}t  \notag \\
&&+\left\langle G\left( \breve{x}_{T}^{i}-\breve{x}_{T}^{(N)}\right) ,\breve{%
x}_{T}^{i}-\breve{x}_{T}^{(N)}\right\rangle -\left\langle G\left( \alpha
_{T}^{i}-\mathbb{E}\alpha _{T}^{i}\right) ,\alpha _{T}^{i}-\mathbb{E}\alpha
_{T}^{i}\right\rangle \Bigg ].  \label{ee7}
\end{eqnarray}%
%
%
%
%
%
%
%
%
%
%
%
%
%
%
%
%
%
%
%
%
%
%
%
%
%
%
%
%
%
%
%
%
%
%
%
%
%
%
%
%
%
%
%
We will use the following
\begin{equation*}
\begin{aligned} &\langle Q(a-b), a-b\rangle-\langle Q(c-d), c-d\rangle\\
=&\langle Q(a-b-(c-d)), a-b-(c-d)\rangle+2\langle Q(a-b-(c-d)), c-d\rangle,
\end{aligned}
\end{equation*}%
and Lemma \ref{lemma for xN and m}, Lemma \ref{lemma for x and xi} as well
as $\mathbb{E}\sup_{0\leq t\leq T}\left\vert \alpha ^{i}(t)\right\vert
^{2}\leq C_{0}$, for some constant $C_{0}$ independent of $N$ which may vary
line by line in the following, we have%
\begin{eqnarray*}
&&\Bigg |\mathbb{E}\Bigg [\int_{0}^{T}\Big [\left\langle Q_{t}\left( \breve{x%
}_{t}^{i}-\breve{x}_{t}^{(N)}\right) ,\breve{x}_{t}^{i}-\breve{x}%
_{t}^{(N)}\right\rangle -\left\langle Q_{t}\left( \alpha _{t}^{i}-\mathbb{E}%
\alpha _{t}^{i}\right) ,\alpha _{t}^{i}-\mathbb{E}\alpha
_{t}^{i}\right\rangle \mathrm{d}t\Bigg | \\
&\leq &C_{0}\int_{0}^{T}\mathbb{E}\left\vert \breve{x}_{t}^{i}-\breve{x}%
_{t}^{(N)}-\left( \alpha _{t}^{i}-\mathbb{E}\alpha _{t}^{i}\right)
\right\vert ^{2}\mathrm{d}t \\
&&+C_{0}\int_{0}^{T}\mathbb{E}\left\vert \breve{x}_{t}^{i}-\breve{x}%
_{t}^{(N)}-\left( \alpha _{t}^{i}-\mathbb{E}\alpha _{t}^{i}\right)
\right\vert \cdot \left\vert \left( \alpha _{t}^{i}-\mathbb{E}\alpha
_{t}^{i}\right) \right\vert \mathrm{d}t \\
&\leq &C_{0}\int_{0}^{T}\mathbb{E}\left\vert \breve{x}_{t}^{i}-\alpha
_{t}^{i}\right\vert ^{2}\mathrm{d}t+C_{0}\int_{0}^{T}\mathbb{E}\left\vert
\breve{x}_{t}^{(N)}-\mathbb{E}\alpha _{t}^{i}\right\vert ^{2}\mathrm{d}t \\
&\leq &O\left( \frac{1}{\sqrt{N}}\right) .
\end{eqnarray*}%
With similar argument, using (\ref{bsde1}) and (\ref{bsde2}), one can show
that%
\begin{eqnarray*}
\Bigg |\mathbb{E}\left[ \left\langle G\left( \breve{x}_{T}^{i}-\breve{x}%
_{T}^{(N)}\right) ,\breve{x}_{T}^{i}-\breve{x}_{T}^{(N)}\right\rangle
-\left\langle G\left( \alpha _{T}^{i}-\mathbb{E}\alpha _{T}^{i}\right)
,\alpha _{T}^{i}-\mathbb{E}\alpha _{T}^{i}\right\rangle \right] \Bigg | %
&\leq &O\left( \frac{1}{\sqrt{N}}\right) , \\
\Bigg |\mathbb{E}\int_{0}^{T}\left[ \left\langle L_{t}\left( \breve{y}%
_{t}^{i}-\breve{y}_{t}^{(N)}\right) ,\breve{y}_{t}^{i}-\breve{y}%
_{t}^{(N)}\right\rangle -\left\langle L_{t}\left( \theta _{t}^{i}-\mathbb{E}%
\theta _{t}^{i}\right) ,\theta _{t}^{i}-\mathbb{E}\theta
_{t}^{i}\right\rangle \right] \mathrm{d}t\Bigg | &\leq &O\left( \frac{1}{%
\sqrt{N}}\right) .
\end{eqnarray*}%
The proof is completed by noticing (\ref{ee7}). \hfill $\Box $

We will prove the control strategies set $(\bar{u}^{1},\bar{u}^{2},\ldots ,%
\bar{u}^{N})$ is an $\epsilon $-Nash equilibrium for \textbf{Problem (CC)}.
For any fixed $i$, $1\leq i\leq N$, we consider the perturbation control $%
u^{i}\in \mathcal{U}_{ad}^{d,i}$ and we have the following state dynamics ($%
j\neq i$):
\begin{equation}
\left\{
\begin{array}{l}
\mathrm{d}\tilde{x}^{i}=\left[ A\tilde{x}^{i}+Bu^{i}+F\tilde{x}^{(N)}+b%
\right] \mathrm{d}t+\left[ Du^{i}+\sigma \right] \mathrm{d}W_{t}^{i}, \\
\mathrm{d}\tilde{y}^{i}=-\left[ M\tilde{x}^{i}+U\tilde{y}^{i}+H\tilde{x}%
^{(N)}+V\tilde{y}^{(N)}+Ku^{i}+f\right] \mathrm{d}t+\tilde{z}^{i}\mathrm{d}%
W_{t}^{i}, \\
\mathrm{d}\tilde{x}^{j}=\left[ A\tilde{x}^{j}+B\varphi \left( \chi
^{j},\beta ^{j},\gamma ^{j}\right) +F\tilde{x}^{(N)}+b\right] \mathrm{d}t+%
\left[ D\varphi \left( \chi ^{j},\beta ^{j},\gamma ^{j}\right) +\sigma %
\right] \mathrm{d}W_{t}^{i}, \\
\mathrm{d}\tilde{y}^{j}=-\left[ M\tilde{x}^{j}+U\tilde{y}^{j}+H\tilde{x}%
^{(N)}+V\tilde{y}^{(N)}+K\varphi \left( \chi ^{j},\beta ^{j},\gamma
^{j}\right) +f\right] \mathrm{d}t+\tilde{z}^{j}\mathrm{d}W_{t}^{j}, \\
\mathrm{d}\alpha ^{j}=\left[ A\alpha ^{j}+B\varphi \left( \chi ^{j},\beta
^{j},\gamma ^{j}\right) +F\mathbb{E}\alpha ^{j}+b\right] \mathrm{d}t+\left[
D\varphi \left( \chi ^{j},\beta ^{j},\gamma ^{j}\right) +\sigma \right]
\mathrm{d}W_{t}^{j}, \\
\mathrm{d}\theta ^{j}=-\left[ M\alpha ^{j}+U\theta ^{j}+H\mathbb{E}\alpha
^{j}+V\mathbb{E}\theta ^{j}+K\varphi \left( \chi ^{j},\beta ^{j},\gamma
^{j}\right) +f\right] \mathrm{d}t+\kappa ^{i}\mathrm{d}W_{t}^{j}, \\
\tilde{x}_{0}^{j}=\tilde{x}_{0}^{i}=\alpha _{0}^{j}=x,\text{ }\tilde{y}%
_{T}^{i}=\Phi \tilde{x}_{T}^{i},\text{ }\tilde{y}_{T}^{j}=\Phi \tilde{x}%
_{T}^{j},\text{ }\theta _{T}^{j}=\Phi \alpha _{T}^{j},\text{ }%
\end{array}%
\right.  \label{ee9}
\end{equation}%
where $\tilde{x}^{(N)}=\displaystyle\frac{1}{N}\sum_{i=1}^{N}\tilde{x}^{i},$
$\tilde{y}^{(N)}=\displaystyle\frac{1}{N}\sum_{i=1}^{N}\tilde{y}^{i}$. The
wellposedness of above system is easily to obtain. To prove $(\bar{u}^{1},%
\bar{u}^{2},\ldots ,\bar{u}^{N})$ is an $\epsilon $-Nash equilibrium, we can
show that for $1\leq i\leq N$,
\begin{equation*}
\inf_{u^{i}\in \mathcal{U}_{ad}^{d,i}}\mathcal{J}_{i}(u^{i},\bar{u}%
_{-i})\geq \mathcal{J}_{i}(\bar{u}^{i},\bar{u}_{-i})-\epsilon .
\end{equation*}%
Then we only need to consider the perturbation $u^{i}\in \mathcal{U}%
_{ad}^{d,i}$ such that $\mathcal{J}_{i}(u^{i},\bar{u}_{-i})\leq \mathcal{J}%
_{i}(\bar{u}^{i},\bar{u}_{-i})$. Thus we have
\begin{equation*}
\mathbb{E}\left[ \int_{0}^{T}\langle Ru^{i}(t),u^{i}(t)\rangle dt\right]
\leq \mathcal{J}_{i}(u^{i},\bar{u}_{-i})\leq \mathcal{J}_{i}(\bar{u}^{i},%
\bar{u}_{-i})\leq J_{i}(\bar{u}^{i})+O\Big(\frac{1}{\sqrt{N}}\Big),
\end{equation*}%
which implies that
\begin{equation}
\mathbb{E}\int_{0}^{T}|u^{i}(t)|^{2}dt\leq C_{0},
\label{boundness of control}
\end{equation}%
where $C_{0}$ is a constant independent of $N$.

Now, for the $i^{th}$ agent, we consider the perturbation in the \textbf{%
Problem (LCC)}. We introduce the following system of the decentralized
limiting state with perturbation control ($j\neq i$):%
\begin{equation}
\left\{
\begin{array}{l}
\mathrm{d}\mathring{x}^{i}=\left[ A\mathring{x}^{i}+Bu^{i}+F\mathbb{E}\alpha
^{i}+b\right] \mathrm{d}t+\left[ Du^{i}+\sigma \right] \mathrm{d}W_{t}^{i},
\\
\mathrm{d}\mathring{y}^{i}=-\left[ M\mathring{x}^{i}+U\mathring{y}^{i}+H%
\mathbb{E}\alpha ^{i}+V\mathbb{E}\mathring{\theta}^{i}+Ku^{i}+f\right]
\mathrm{d}t+\mathring{z}^{i}\mathrm{d}W_{t}^{i}, \\
\mathrm{d}\alpha ^{i}=\left[ A\alpha ^{i}+B\varphi \left( \chi ^{i},\beta
^{i},\gamma ^{i}\right) +F\mathbb{E}\alpha ^{i}+b\right] \mathrm{d}t+\left[
D\varphi \left( \chi ^{i},\beta ^{i},\gamma ^{i}\right) +\sigma \right]
\mathrm{d}W_{t}^{i}, \\
\mathrm{d}\theta ^{i}=-\left[ M\alpha ^{i}+U\theta ^{i}+H\mathbb{E}\alpha
^{i}+V\mathbb{E}\theta ^{i}+K\varphi \left( \chi ^{i},\beta ^{i},\gamma
^{i}\right) +f\right] \mathrm{d}t+\kappa ^{i}\mathrm{d}W_{t}^{i}, \\
\mathring{x}_{0}^{i}=\alpha _{0}^{i}=x,\text{ }\mathring{y}_{T}^{i}=\Phi
\mathring{x}_{T}^{i},\text{ }\theta _{T}^{i}=\Phi \alpha _{T}^{i}.%
\end{array}%
\right.  \label{ee10}
\end{equation}%
We have the following results:

\begin{lemma}
\label{lemma 2}Let \emph{(A1)} and \emph{(A2)} hold, then
\begin{eqnarray}
\mathbb{E}\left[ \sup_{0\leq t\leq T}\Big|\tilde{x}_{t}^{(N)}-\mathbb{E}%
\alpha ^{i}(t)\Big|^{2}\right] &=&O\Big(\frac{1}{N}\Big),  \label{perineq1}
\\
\mathbb{E}\left[ \sup_{0\leq t\leq T}\Big|\tilde{y}_{t}^{(N)}-\mathbb{E}%
\theta ^{i}(t)\Big|^{2}\right] &=&O\Big(\frac{1}{N}\Big).  \label{perineq2}
\end{eqnarray}
\end{lemma}

\paragraph{Proof.}

By (\ref{ee9}), we get
\begin{equation}
\left\{
\begin{array}{l}
\mathrm{d}\tilde{x}^{\left( N\right) }=\Big [\left( A+F\right) \tilde{x}%
^{(N)}+\frac{1}{N}Bu^{i}+\frac{1}{N}\sum_{j=1,j\neq i}^{N}B\varphi \left(
\chi ^{j},\beta ^{j},\gamma ^{j}\right) \\
\qquad \qquad +b\Big ]\mathrm{d}t+\frac{1}{N}Du^{i}W_{t}^{i}+\frac{1}{N}%
\sum_{j=1}^{N}\sigma \mathrm{d}W_{t}^{j} \\
\qquad \qquad +\frac{1}{N}\sum_{j=1,j\neq i}^{N}B\varphi \left( \chi
^{j},\beta ^{j},\gamma ^{j}\right) \mathrm{d}W_{t}^{j}, \\
\mathrm{d}\tilde{y}^{\left( N\right) }=-\Big [\left( M+H\right) \tilde{x}%
^{(N)}+\left( U+V\right) \tilde{y}^{(N)}+\frac{1}{N}Ku^{i} \\
\qquad \qquad +\frac{1}{N}\sum_{j=1,j\neq i}^{N}B\varphi \left( \chi
^{j},\beta ^{j},\gamma ^{j}\right) +f\Big ]\mathrm{d}t+\frac{1}{N}%
\sum_{j=1}^{N}\tilde{z}^{j}\mathrm{d}W_{t}^{j}, \\
\tilde{x}_{0}^{\left( N\right) }=x,\text{ }\tilde{y}_{T}^{\left( N\right)
}=\Phi \tilde{x}_{T}^{\left( N\right) }.%
\end{array}%
\right.  \label{ee12}
\end{equation}%
Let us denote
\begin{eqnarray*}
\Pi &:&=\tilde{x}^{\left( N\right) }-\mathbb{E}\alpha ^{i}, \\
\Lambda &:&=\tilde{y}^{\left( N\right) }-\mathbb{E}\theta ^{i},
\end{eqnarray*}%
and recall (\ref{ee2}) which is%
\begin{equation*}
\left\{
\begin{array}{l}
\mathrm{d}\mathbb{E}\alpha ^{i}=\left[ \left( A+F\right) \mathbb{E}\alpha
^{i}+B\mathbb{E}\varphi \left( \chi ^{i},\beta ^{i},\gamma ^{i}\right) +b%
\right] \mathrm{d}t, \\
\mathrm{d}\mathbb{E}\theta ^{i}=-\left[ \left( M+H\right) \mathbb{E}\alpha
^{i}+\left( U+V\right) \mathbb{E}\theta ^{i}+K\mathbb{E}\varphi \left( \chi
^{i},\beta ^{i},\gamma ^{i}\right) +f\right] \mathrm{d}t, \\
\mathbb{E}\alpha _{0}^{i}=x,\text{ }\mathbb{E}\theta _{T}^{i}=\Phi \mathbb{E}%
\alpha _{T}^{i},\text{ }%
\end{array}%
\right.
\end{equation*}%
we have
\begin{equation*}
\left\{
\begin{array}{l}
\mathrm{d}\Pi =\Big [\left( A+F\right) \Pi +\frac{1}{N}Bu^{i} \\
\qquad \qquad +\left( \frac{1}{N}\sum_{j=1,i\neq j}^{N}B\varphi \left( \chi
^{i},\beta ^{i},\gamma ^{i}\right) -\mathbb{E}B\varphi \left( \chi
^{i},\beta ^{i},\gamma ^{i}\right) \right) \Big ]\mathrm{d}t \\
\qquad \qquad +\frac{1}{N}Du_{i}W_{t}^{i}+\frac{1}{N}\sum_{j=1}^{N}\sigma
\mathrm{d}W_{t}^{j} \\
\qquad \qquad +\frac{1}{N}\sum_{j=1,j\neq i}^{N}B\varphi \left( \chi
^{j},\beta ^{j},\gamma ^{j}\right) \mathrm{d}W_{t}^{j}, \\
\mathrm{d}\Lambda =\Big [\left( M+H\right) \Pi +\left( U+V\right) \Lambda +%
\frac{1}{N}Ku^{i} \\
\qquad \qquad +\left( \frac{1}{N}\sum_{j=1,i\neq j}^{N}B\varphi \left( \chi
^{i},\beta ^{i},\gamma ^{i}\right) -\mathbb{E}B\varphi \left( \chi
^{i},\beta ^{i},\gamma ^{i}\right) \right) \Big ]\mathrm{d}t \\
\qquad \qquad +\frac{1}{N}\sum_{j=1}^{N}\tilde{z}^{i}\mathrm{d}W_{t}^{i} \\
\Pi _{0}=0,\text{ }\Lambda _{T}=\Phi \Pi _{T}.%
\end{array}%
\right.
\end{equation*}%
By the Cauchy-Schwartz inequality as well as the B-D-G inequality, we obtain
that there exists a constant $C_{0}$ independent of $N$ which may vary line
by line such that, for any $t\in \lbrack 0,T]$,
\begin{equation}
\begin{aligned} \mathbb{E}\sup_{0\leq s\leq t}|\Pi_s|^2\leq & \,
C_0\mathbb{E}\int_0^t\left(|\Pi_s|^2+\frac{1}{N^2}|u^i_s|^2\right)ds \\
&+C_0\mathbb{E}\int_0^t\left|\frac{1}{N}\sum\limits_{j=1,j\neq
i}^{N}\varphi(\chi ^{j},\beta ^{j},\gamma ^{j}) -\mathbb{E}\varphi(\chi
^{j},\beta ^{j},\gamma ^{j})\right|^2ds \\
&+\frac{C_0}{N^2}\mathbb{E}\sum_{j=1}^{N}\int_0^t\left|\sigma_s\right|^2ds
\\
&+\frac{C_0}{N^2}\mathbb{E}\int_0^t|u^i_s|^2ds+\frac{C_0}{N^2}\mathbb{E}\!\!%
\!\sum\limits_{j=1,j\neq i}^{N}\int_0^t|\varphi(\chi ^{j},\beta ^{j},\gamma
^{j}|^2ds. \end{aligned}  \label{ee11}
\end{equation}%
On the one hand, by denoting $\mu (s):=\mathbb{E}\varphi (\chi ^{j},\beta
^{j},\gamma ^{j})$ (note that since $(\chi ^{j},\beta ^{j},\gamma ^{j}$, $%
1\leq j\leq N$, $j\neq i$, are independent identically distributed, thus $%
\mu $ is independent of $j$), we have
\begin{equation*}
\begin{aligned} &\mathbb{E}\left|\frac{1}{N}\sum\limits_{j=1,j\neq
i}^{N}\varphi(\chi ^{j},\beta ^{j},\gamma ^{j}) -\mu_s\right|^2\\ \leq& \,
2\mathbb{E}\left|\frac{1}{N}\sum\limits_{j=1,j\neq i}^{N}\varphi(\chi
^{j},\beta ^{j},\gamma ^{j}) -\frac{N-1}{N}\mu_s\right|^2
+2\mathbb{E}\left|\frac{1}{N}\mu_s\right|^2\\ =& \,
2\frac{(N-1)^2}{N^2}\mathbb{E}\left|\frac{1}{N-1}\sum\limits_{j=1,j\neq
i}^{N}\varphi(\chi ^{j},\beta ^{j},\gamma ^{j})
-\mu_s\right|^2+\frac{2}{N^2}\mathbb{E}|\mu_s|^2. \end{aligned}
\end{equation*}%
Then, due to the fact that $(\chi ^{i},\beta ^{i},\gamma ^{i}),1\leq i\leq N$
are identically distributed and $\varphi (\chi ^{i},\beta ^{i},\gamma
^{i})\in L_{\mathbb{F}^{i}}^{2}(0,T;U)$, similarly to Lemma \ref{lemma for
xN and m} we can obtain that there exists a constant $C_{0}$ independent of $%
N$ such that
\begin{equation*}
\begin{aligned} &\int_0^t\mathbb{E}\left|\frac{1}{N}\sum\limits_{j=1,j\neq
i}^{N}\varphi(\chi ^{j},\beta ^{j},\gamma ^{j}) -\mathbb{E}\varphi(\chi
^{j},\beta ^{j},\gamma ^{j})\right|^2 \\ \leq & \,
\frac{C_0(N-1)^2}{N^2}\int_0^t\mathbb{E}\left|\frac{1}{N-1}\sum%
\limits_{j=1,j\neq i}^{N}\varphi(\chi ^{j},\beta ^{j},\gamma
^{j})-\mu_s\right|^2ds +\frac{C_0}{N^2}\int_0^t\mathbb{E}|\mu_s|^2ds \\ =&
\, \frac{C_0(N-1)}{N^2}\int_0^t\mathbb{E}\left|\varphi(\chi ^{j},\beta
^{j},\gamma ^{j})
-\mu_s\right|^2ds+\frac{C_0}{N^2}\int_0^t\mathbb{E}|\mu_s|^2ds \\ =& \,
O\Big(\frac{1}{N}\Big). \end{aligned}
\end{equation*}%
In addition, due to (\ref{boundness of control}), we get
\begin{equation*}
\begin{aligned} \frac{C_0}{N^2}\mathbb{E}\int_0^t|u^i_s|^2ds
+\frac{C_0}{N^2}\mathbb{E}\sum_{j=1}^{N}\int_0^t\left|\sigma_s\right|^2ds =
O\Big(\frac{1}{N}\Big), \end{aligned}
\end{equation*}%
and similarly, since $(\chi ^{j},\beta ^{j},\gamma ^{j})$, $1\leq j\leq N$, $%
j\neq i$, are identically distributed, we have
\begin{equation*}
\frac{C_{0}}{N^{2}}\mathbb{E}\sum\limits_{j=1,j\neq
i}^{N}\int_{0}^{t}|\varphi (\chi ^{j},\beta ^{j},\gamma ^{j}|^{2}ds=O\Big(%
\frac{1}{N}\Big).
\end{equation*}%
Therefore, from above estimates, we get from (\ref{ee11}) that, for any $%
t\in \lbrack 0,T]$,
\begin{equation*}
\mathbb{E}\sup_{0\leq s\leq t}|\Pi _{s}|^{2}\leq C_{0}\mathbb{E}%
\int_{0}^{t}|\Pi _{s}|^{2}ds+O\Big(\frac{1}{N}\Big).
\end{equation*}%
Finally, by using Gronwall's inequality, we get (\ref{perineq1}). We now
proceed the second inequality. Applying It\^{o}'s formula again, we have%
\begin{eqnarray*}
\left\vert \Lambda _{t}\right\vert ^{2}+\frac{1}{N}\sum_{j=1}^{N}\tilde{z}%
_{t}^{i}\mathrm{d}t &=&\left\vert \Lambda _{T}\right\vert ^{2}+2\int_{t}^{T}%
\Big <\Lambda _{t},\left( M+H\right) \Pi +\left( U+V\right) \Lambda +\frac{1%
}{N}Ku^{i} \\
&&+\left( \frac{1}{N}\sum_{j=1,i\neq j}^{N}B\varphi \left( \chi ^{i},\beta
^{i},\gamma ^{i}\right) -\mathbb{E}B\varphi \left( \chi ^{i},\beta
^{i},\gamma ^{i}\right) \right) \Big >\mathrm{d}s \\
&&-2\frac{1}{N}\sum_{i=1}^{N}\int_{t}^{T}\Big <\Lambda _{t},\breve{z}^{i}%
\Big >\mathrm{d}W_{s}^{i}.
\end{eqnarray*}%
Still using B-D-G inequalities, we show that there exists a constant $C_{2},$
modifying $C_{2}$ if necessary,%
\begin{eqnarray*}
&&\mathbb{E}\left[ \sup_{0\leq t\leq T}\left \vert \Lambda _{t}\right \vert
^{2}+\frac{1}{N}\sum_{i=1}^{N}\int_{0}^{T}\left \vert \breve{z}^{i}\right
\vert ^{2}\mathrm{d}s\right] \\
&\leq &\mathbb{E}\left \vert \Phi \Lambda _{T}\right \vert ^{2}+C_{2}\mathbb{%
E}\Bigg [\int_{0}^{T}\Big <\Lambda _{t},\left( M+H\right) \Pi \Big > \\
&&+\Big <\Lambda _{t},\left( U+V\right) \Lambda \Big >+\Big <\Lambda _{t},%
\frac{1}{N}Ku^{i}\Big > \\
&&+\Big <\Lambda _{t},\frac{1}{N}\sum_{i=1}^{N}K\varphi \left( \chi
^{i},\beta ^{i},\gamma ^{i}\right) -K\mathbb{E}\varphi \left( \chi
^{i},\beta ^{i},\gamma ^{i}\right) \Big >\mathrm{d}s \\
&&+\frac{1}{N}\sum_{i=1}^{N}\left( \int_{0}^{T}\left \vert \Lambda
_{t}\right \vert ^{2}\left \vert \breve{z}_{t}^{i}\right \vert ^{2}\mathrm{d}%
t\right) ^{\frac{1}{2}}\Bigg ] \\
&\leq &\mathbb{E}\left \vert \Phi \Delta _{T}^{1}\right \vert ^{2}+C_{1}%
\mathbb{E}\Bigg [\int_{t}^{T}\left \vert \Lambda _{t}\right \vert \Big (%
\left \vert \left( M+H\right) \Pi \right \vert +\left \vert \left(
U+V\right) \Lambda \right \vert \\
&&+\left \vert \frac{1}{N}Ku^{i}\right \vert +\left \vert \frac{1}{N}%
\sum_{i=1}^{N}K\varphi \left( \chi ^{i},\beta ^{i},\gamma ^{i}\right) -K%
\mathbb{E}\varphi \left( \chi ^{i},\beta ^{i},\gamma ^{i}\right) \right
\vert \Big )\mathrm{d}t \\
&&+\frac{1}{N}\sum_{i=1}^{N}\left( \int_{0}^{T}\left \vert \Lambda
_{t}\right \vert ^{2}\left \vert \breve{z}^{i}\right \vert ^{2}\mathrm{d}%
t\right) ^{\frac{1}{2}}\Bigg ].
\end{eqnarray*}

By employing the classical Cauchy-Schwarz inequality and Gronwall's
inequality, we obtain get (\ref{perineq2}).\hfill $\Box $

\begin{lemma}
\label{lemma 3}
\begin{eqnarray}
\mathbb{E}\left[ \sup_{0\leq t\leq T}\Big|\tilde{x}_{t}^{i}-\mathring{x}%
_{t}^{i}\Big|^{2}\right] &=&O\Big(\frac{1}{N}\Big),  \label{per1} \\
\mathbb{E}\left[ \sup_{0\leq t\leq T}\Big|\tilde{y}_{t}^{i}-\mathring{y}%
_{t}^{i}\Big|^{2}\right] &=&O\Big(\frac{1}{N}\Big).  \label{per2}
\end{eqnarray}
\end{lemma}

\paragraph{Proof.}

From respectively the first equation of (\ref{ee9}) and (\ref{ee10}), we
obtain%
\begin{equation*}
\left\{
\begin{array}{l}
\mathrm{d}\left( \tilde{x}^{i}-\mathring{x}^{i}\right) =\left[ A\left(
\tilde{x}^{i}-\mathring{x}^{i}\right) +F\left( \tilde{x}^{(N)}-\mathbb{E}%
\alpha ^{i}\right) \right] \mathrm{d}t \\
\mathrm{d}\left( \tilde{y}^{i}-\mathring{y}^{i}\right) =-\Big [M\left(
\tilde{x}^{i}-\mathring{x}^{i}\right) +U\left( \tilde{y}^{i}-\mathring{y}%
^{i}\right) \\
\qquad \qquad +H\left( \tilde{x}^{(N)}-\mathbb{E}\alpha ^{i}\right) +V\left(
\tilde{y}^{(N)}-\mathbb{E}\theta ^{i}\right) \Big ]\mathrm{d}t+\left( \tilde{%
z}^{i}-\mathring{z}^{i}\right) \mathrm{d}W_{t}^{i}, \\
\tilde{x}_{0}^{i}-\mathring{x}_{0}^{i}=0,\text{ }\tilde{y}_{T}^{i}-\mathring{%
y}_{T}^{i}=\Phi \left( \tilde{x}_{T}^{i}-\mathring{x}_{T}^{i}\right) .%
\end{array}%
\right.
\end{equation*}
With the help of classical estimates of SDE and BSDE, Gronwall's inequality
and (\ref{perineq1}) and (\ref{perineq2}) of Lemma \ref{lemma 2}, it is
easily to obtain (\ref{per1}) and (\ref{per2}). The proof is completed.
\hfill $\Box $

\begin{lemma}
\label{lemma 4} For all $1\leq i\leq N,$ for the perturbation control $u^{i}$%
, we have
\begin{equation*}
\Big|\mathcal{J}_{i}({u}^{i},\bar{u}_{-i})-J_{i}({u}^{i})\Big|=O\Big(\frac{1%
}{\sqrt{N}}\Big).
\end{equation*}
\end{lemma}

\paragraph{Proof.\protect\medskip}

Recall (\ref{cost1}), (\ref{cost2}), (\ref{consi1}), and (\ref{consi2}), we
have%
\begin{eqnarray*}
&&\mathcal{J}_{i}\left( u^{i},u_{-i}\right) -J_{i}\left( u^{i}\right) \\
&=&\frac{1}{2}\mathbb{E}\Bigg [\int_{0}^{T}\Big [\left\langle Q_{t}\left(
\tilde{x}_{t}^{i}-\tilde{x}_{t}^{(N)}\right) ,\tilde{x}_{t}^{i}-\tilde{x}%
^{(N)}\right\rangle +\left\langle L_{t}\left( \tilde{y}_{t}^{i}-\tilde{y}%
_{t}^{(N)}\right) ,\tilde{y}_{t}^{i}-\tilde{y}_{t}^{(N)}\right\rangle \\
&&-\left\langle Q_{t}\left( \mathring{x}_{t}^{i}-\mathbb{E}\alpha
_{t}^{i}\right) ,\mathring{x}_{t}^{i}-\mathbb{E}\alpha _{t}^{i}\right\rangle
-\left\langle L_{t}\left( \mathring{y}_{t}^{i}-\mathbb{E}\theta
_{t}^{i}\right) ,\mathring{y}_{t}^{i}-\mathbb{E}\theta _{t}^{i}\right\rangle %
\Big ]\mathrm{d}t \\
&&+\left\langle G\left( x_{T}^{i}-\phi _{T}^{1}\right) ,x_{T}^{i}-\phi
_{T}^{1}\right\rangle -\left\langle G\left( \mathring{x}_{T}^{i}-\mathbb{E}%
\alpha _{T}^{i}\right) ,x_{T}^{i}-\mathbb{E}\alpha _{T}^{i}\right\rangle
\end{eqnarray*}%
Using Lemma \ref{lemma 2} and Lemma \ref{lemma 3} as well as $\mathbb{E}%
\sup_{0\leq t\leq T}\left( |\bar{y}^{i}(t)|^{2}+|\alpha ^{i}(t)|^{2}\right)
\leq C_{0}$, for some constant $C$ independent of $N$ which may vary line by
line in the following, we have
\begin{eqnarray*}
&&\Bigg |\mathbb{E}\Bigg [\int_{0}^{T}\Big [\left\langle Q_{t}\left( \tilde{x%
}_{t}^{i}-\tilde{x}_{t}^{(N)}\right) ,\tilde{x}_{t}^{i}-\tilde{x}%
^{(N)}\right\rangle +\left\langle L_{t}\left( \tilde{y}_{t}^{i}-\tilde{y}%
_{t}^{(N)}\right) ,\tilde{y}_{t}^{i}-\tilde{y}_{t}^{(N)}\right\rangle \\
&&-\left\langle Q_{t}\left( \mathring{x}_{t}^{i}-\mathbb{E}\alpha
_{t}^{i}\right) ,\mathring{x}_{t}^{i}-\mathbb{E}\alpha _{t}^{i}\right\rangle
-\left\langle L_{t}\left( \mathring{y}_{t}^{i}-\mathbb{E}\theta
_{t}^{i}\right) ,\mathring{y}_{t}^{i}-\mathbb{E}\theta _{t}^{i}\right\rangle %
\Big ]\mathrm{d}t\Bigg | \\
&\leq &C\int_{0}^{T}\mathbb{E}\left\vert \tilde{x}_{t}^{i}-\tilde{x}%
_{t}^{(N)}-\left( \mathring{x}_{t}^{i}-\mathbb{E}\alpha _{t}^{i}\right)
\right\vert ^{2}\mathrm{d}t \\
&&+C\int_{0}^{T}\mathbb{E}\left\vert \tilde{x}_{t}^{i}-\tilde{x}%
_{t}^{(N)}-\left( \mathring{x}_{t}^{i}-\mathbb{E}\alpha _{t}^{i}\right)
\right\vert \cdot \left\vert \mathring{x}_{t}^{i}-\mathbb{E}\alpha
_{t}^{i}\right\vert \mathrm{d}t \\
&&+C\int_{0}^{T}\mathbb{E}\left\vert \tilde{y}_{t}^{i}-\tilde{y}%
_{t}^{(N)}-\left( \mathring{y}_{t}^{i}-\mathbb{E}\theta _{t}^{i}\right)
\right\vert ^{2}\mathrm{d}t \\
&&+C\int_{0}^{T}\mathbb{E}\left\vert \tilde{y}_{t}^{i}-\tilde{y}%
_{t}^{(N)}-\left( \mathring{y}_{t}^{i}-\mathbb{E}\theta _{t}^{i}\right)
\right\vert \cdot \left\vert \mathring{y}_{t}^{i}-\mathbb{E}\theta
_{t}^{i}\right\vert \mathrm{d}t \\
&\leq &C\int_{0}^{T}\mathbb{E}\left[ \left\vert \tilde{x}_{t}^{i}-\tilde{x}%
_{t}^{(N)}\right\vert ^{2}+\left\vert \mathring{x}_{t}^{i}-\mathbb{E}\alpha
_{t}^{i}\right\vert ^{2}\right] \mathrm{d}t \\
&&+C\int_{0}^{T}\mathbb{E}\left( \left\vert \tilde{x}_{t}^{i}-\tilde{x}%
_{t}^{(N)}\right\vert ^{2}+\left\vert \mathring{x}_{t}^{i}-\mathbb{E}\alpha
_{t}^{i}\right\vert ^{2}\right) ^{\frac{1}{2}} \\
&&+C\int_{0}^{T}\mathbb{E}\left[ \left\vert \tilde{y}_{t}^{i}-\tilde{y}%
_{t}^{(N)}\right\vert ^{2}+\left\vert \mathring{y}_{t}^{i}-\mathbb{E}\theta
_{t}^{i}\right\vert ^{2}\right] \mathrm{d}t
\end{eqnarray*}%
\begin{eqnarray*}
&&+C\int_{0}^{T}\mathbb{E}\left( \left\vert \tilde{y}_{t}^{i}-\tilde{y}%
_{t}^{(N)}\right\vert ^{2}+\left\vert \mathring{y}_{t}^{i}-\mathbb{E}\theta
_{t}^{i}\right\vert ^{2}\right) ^{\frac{1}{2}} \\
&\leq &O\left( \frac{1}{N}\right) .
\end{eqnarray*}%
With similar argument, we can show that
\begin{equation*}
\left\vert \mathbb{E}\left[ \left\langle G\left( x_{T}^{i}-\phi
_{T}^{1}\right) ,x_{T}^{i}-\phi _{T}^{1}\right\rangle -\left\langle G\left(
\mathring{x}_{T}^{i}-\mathbb{E}\alpha _{T}^{i}\right) ,x_{T}^{i}-\mathbb{E}%
\alpha _{T}^{i}\right\rangle \right] \right\vert \leq O\left( \frac{1}{N}%
\right) .
\end{equation*}%
Hence, we get the desired result. \hfill $\Box $


\noindent \textbf{Proof of Theorem \ref{Nash equilibrium theorem}:} Now, we
consider the $\epsilon $-Nash equilibrium for $\mathcal{A}_{i}$ for \textbf{%
Problem (CC)}. Combining Lemma \ref{first lemma for cost} and Lemma \ref%
{lemma 4}, we have
\begin{equation*}
\mathcal{J}_{i}(\bar{u}^{i},\bar{u}_{-i})=J_{i}(\bar{u}^{i})+O\Big(\frac{1}{%
\sqrt{N}}\Big)\leq J_{i}({u}^{i})+O\Big(\frac{1}{\sqrt{N}}\Big)=\mathcal{J}%
_{i}({u}^{i},\bar{u}_{-i})+O\Big(\frac{1}{\sqrt{N}}\Big).
\end{equation*}%
Consequently, Theorem \ref{Nash equilibrium theorem} holds with $\epsilon =O%
\Big(\frac{1}{\sqrt{N}}\Big)$. \hfill $\Box $

\appendix

\section{Proof of theorem}

\label{app1}

\paragraph{Proof of Theorem \protect\ref{th1}.}

(\textbf{Uniqueness}) Suppose that there exists two solutions: $%
(x^{1},y^{1},z^{1},p^{1},q^{1},k^{1}),$

\noindent $(x^{2},y^{2},z^{2},p^{2},q^{2},k^{2})$ and denote%
\begin{eqnarray*}
\hat{x} &=&x^{1}-x^{2},\text{ }\hat{z}=z^{1}-z^{2}, \\
\hat{y} &=&y^{1}-y^{2},\text{ }\hat{p}=p^{1}-p^{2}, \\
\hat{q} &=&q^{1}-q^{2},\text{ }\hat{k}=k^{1}-k^{2},
\end{eqnarray*}%
Then, we have
\begin{equation}
\left\{
\begin{array}{l}
\mathrm{d}\hat{x}=\left[ A\hat{x}+B\hat{\varphi}\left( p,q,k\right) +F%
\mathbb{E}\hat{x}\right] \mathrm{d}t+\left[ D\hat{\varphi}\left( \hat{p},%
\hat{q},\hat{k}\right) \right] \mathrm{d}W_{t}, \\
\mathrm{d}\hat{y}=-\left[ M\hat{x}+U\hat{y}+H\mathbb{E}\hat{x}+V\mathbb{E}%
\hat{y}+K\hat{\varphi}\left( \hat{p},\hat{q},\hat{k}\right) \right] \mathrm{d%
}t+\hat{z}\mathrm{d}W_{t}, \\
\mathrm{d}\hat{p}=\left[ U\hat{p}-L\left( \hat{y}-\mathbb{E}\hat{y}\right) %
\right] \mathrm{d}t, \\
\mathrm{d}\hat{q}=\left[ -M\hat{p}-A\hat{q}+Q\left( \hat{x}-\mathbb{E}\hat{x}%
\right) \right] \mathrm{d}t+\hat{k}\mathrm{d}W_{t}, \\
\hat{x}_{0}=0\in \mathbb{R}^{n},\text{ }\hat{y}_{T}=\Phi \hat{x}_{T},\text{ }%
\hat{p}_{0}=0, \\
\hat{q}_{T}=\Phi ^{T}\hat{p}_{T}-G\left( x_{T}-\mathbb{E}x_{T}\right) .%
\end{array}%
\right.  \label{e221}
\end{equation}%
with
\begin{eqnarray*}
\hat{\varphi}\left( \hat{p},\hat{q},\hat{k}\right) &=&\varphi \left(
p^{1},q^{1},k^{1}\right) -\varphi \left( p^{2},q^{2},k^{2}\right) \\
&=&\mathbf{P}_{\Gamma }\left[ R^{-1}\left(
B^{T}q^{1}+K^{T}p^{1}+D^{T}k^{1}\right) \right] -\mathbf{P}_{\Gamma }\left[
R^{-1}\left( B^{T}q^{2}+K^{T}p^{2}+D^{T}k^{2}\right) \right]
\end{eqnarray*}%
Taking the expectation in the second equation of (\ref{e221}) yields $%
\mathbb{E}\hat{p}=0$. Applying It\^{o}'s formula to $\big<\hat{q},\hat{x}%
\big>-\big<\hat{p},\hat{y}\big>$ and taking expectations on both sides
(also, noting $\mathbb{E}\hat{p}=0,$ which derives that $\mathbb{E}\hat{q}%
=0, $ and the monotonicity property of $\widehat{\varphi }$)$,$ we arrive at
\begin{eqnarray*}
0 &=&\mathbb{E}\left[ \left\langle G\left( \hat{x}_{T}-\mathbb{E}\hat{x}%
_{T}\right) ,\hat{x}_{T}\right\rangle +\Psi \hat{y}_{0}\left( \hat{y}_{0}-%
\mathbb{E}\hat{y}_{0}\right) \right] \\
&&+\mathbb{E}\int_{0}^{T}\left\langle B^{T}\hat{q}_{s}+D^{T}\hat{k}_{s}+K^{T}%
\hat{p}_{s},\hat{\varphi}\left( \hat{p},\hat{q},\hat{k}\right) \right\rangle
\mathrm{d}s \\
&&+\mathbb{E}\int_{0}^{T}\left[ \left\langle \hat{x}_{s},Q\left( \hat{x}_{s}-%
\mathbb{E}\hat{x}_{s}\right) \right\rangle +L\hat{y}_{s}\cdot \left( \hat{y}%
_{s}-\mathbb{E}\hat{y}_{s}\right) \right] \mathrm{d}s \\
&&+\mathbb{E}\int_{0}^{T}\left[ \left\langle \hat{p}_{s},H\mathbb{E}\hat{x}%
_{s}\right\rangle +\hat{p}_{s}V\mathbb{E}\hat{y}_{s}+\left\langle \hat{q}%
_{s},F\mathbb{E}\hat{x}_{s}\right\rangle \right] \mathrm{d}s \\
&\geq &\mathbb{E}\left[ \left\langle G^{\frac{1}{2}}\left( \hat{x}_{T}-%
\mathbb{E}\hat{x}_{T}\right) ,G^{\frac{1}{2}}\left( \hat{x}_{T}-\mathbb{E}%
\hat{x}_{T}\right) \right\rangle \right] \\
&&+\mathbb{E}\int_{0}^{T}\left[ \left\langle Q^{\frac{1}{2}}\left( \hat{x}%
_{s}-\mathbb{E}\hat{x}_{s}\right) ,Q^{\frac{1}{2}}\left( \hat{x}_{s}-\mathbb{%
E}\hat{x}_{s}\right) \right\rangle +L\left( \hat{y}_{s}-\mathbb{E}\hat{y}%
_{s}\right) \cdot \left( \hat{y}_{s}-\mathbb{E}\hat{y}_{s}\right) \right]
\mathrm{d}s.
\end{eqnarray*}%
Thus, $G\big(\hat{x}_{T}-\mathbb{E}\hat{x}_{T}\big)=0$, $Q\big(\hat{x}-%
\mathbb{E}\hat{x}\big)=0,$ and $L\left( \hat{y}_{s}-\mathbb{E}\hat{y}%
_{s}\right) =0$ which according to the uniqueness and existences of
classical BSDE theory implies $\hat{p}_{s}\equiv 0,$ $\hat{q}_{s}\equiv 0.$
Next, we have $\widehat{\varphi }(\hat{p},\hat{q},\hat{k})\equiv 0$ which
further implies $\mathbb{E}\hat{x}_{s}\equiv 0,$ hence $\hat{x}_{s}\equiv 0.$
Moreover, $\hat{y}_{T}=0$ yields, by Theorem 3.1 in \cite{BLP}, $\hat{y}=0.$
Hence the uniqueness follows. ~\hfill $\Box $

\vskip12pt

In order to prove the existence for FBSDE (\ref{cc}), we need the following
result. It involves a priori estimates of solutions of the following family
of mean fields FBSDEs parameterized by $\alpha \in \lbrack 0,1].$

Before that, we denote
\begin{eqnarray*}
\mathbf{B} &\mathbf{\triangleq }&Ax+B\varphi \left( p,q,k\right) +F\mathbb{E}%
x+b, \\
\Sigma &\mathbf{\triangleq }&D\varphi \left( p,q,k\right) +\sigma , \\
\mathbf{F} &\mathbf{\triangleq }&-\left[ Mx+Uy+H\mathbb{E}x+V\mathbb{E}%
y+K\varphi \left( p,q,k\right) +f\right] , \\
\Xi &\mathbf{\triangleq }&Up-L\left( y-\mathbb{E}y\right) , \\
\Upsilon &\mathbf{\triangleq }&-Mp-Aq+Q\left( x-\mathbb{E}x\right) ,
\end{eqnarray*}%
Consider the following a family of FBSDEs with parameter $\alpha \in \mathbb{%
R},$%
\begin{equation}
\left\{
\begin{array}{l}
\mathrm{d}x=\left[ \alpha \mathbf{B}\left( x,\mathbb{E}x,p,q,k\right) +b_{0}%
\right] \mathrm{d}t+\left[ \alpha \Sigma \left( p,q,k\right) +\sigma _{0}%
\right] \mathrm{d}W_{t}, \\
\mathrm{d}y=\left[ \alpha \mathbf{F}\left( x,\mathbb{E}x,y,\mathbb{E}%
y,p,q,k\right) +\gamma -\mathbb{E\gamma }\right] \mathrm{d}t+z\mathrm{d}%
W_{t}, \\
\mathrm{d}p=\left[ \alpha \Xi \left( y,\mathbb{E}y,p\right) +\lambda _{0}-%
\mathbb{E}\lambda _{0}\right] \mathrm{d}t, \\
\mathrm{d}q=\left[ \alpha \Upsilon \left( x,\mathbb{E}x,p,q\right) +\psi
_{0}-\mathbb{E}\psi _{0}\right] \mathrm{d}t+k\mathrm{d}W_{t}, \\
x_{0}=x\in \mathbb{R}^{n},\text{ }y_{T}=\alpha \Phi x_{T}+\upsilon -\mathbb{%
E\upsilon }, \\
p_{0}=\varsigma -\mathbb{E\varsigma }, \\
q_{T}=\alpha \Phi ^{T}\hat{p}_{T}-\alpha G\left( x_{T}-\mathbb{E}%
x_{T}\right) +\xi -\mathbb{E\xi },%
\end{array}%
\right.  \label{pfbsde}
\end{equation}%
where $\left( b_{0},\sigma _{0},\gamma _{0},\lambda _{0},\mu _{0},\psi
_{0}\right) \in L_{\mathbb{F}^{W}}^{2}(0,T;\mathbb{R}^{n}\times \mathbb{R}%
^{n}\times \mathbb{R\times R}^{n}\mathbb{\times R}^{n}\times \mathbb{R}),$
and $\upsilon $ ($\xi )$ is a $\mathbb{R}$-valued ($\mathbb{R}^{n}$-valued)
square integrable random varible which is $\mathbb{F}_{T}^{W}$-measurable.
Note the coefficient $\Xi \triangleq Up-L\left( y-\mathbb{E}y\right) .$ It
is easy to check $\mathbb{E}p^{\alpha }=0,$ then by uniqueness of BSDE, $%
\mathbb{E}q^{\alpha }=0.$ Specifically, letting $\alpha =0,$ one immediately
has
\begin{equation}
\left\{
\begin{array}{l}
\mathrm{d}x=b_{0}\mathrm{d}t+\sigma _{0}\mathrm{d}W_{t}, \\
\mathrm{d}y=\left( \gamma _{0}-\mathbb{E\gamma }_{0}\right) \mathrm{d}t+z%
\mathrm{d}W_{t}, \\
\mathrm{d}p=\lambda _{0}\mathrm{d}t, \\
\mathrm{d}q=\left( \psi _{0}-\mathbb{E}\psi _{0}\right) \mathrm{d}t+k\mathrm{%
d}W_{t}, \\
x_{0}^{\alpha }=x\in \mathbb{R}^{n},\text{ }y_{T}^{\alpha }=\upsilon -%
\mathbb{E\upsilon }, \\
p_{0}^{\alpha }=\varsigma -\mathbb{E\varsigma },\text{ }q_{T}^{\alpha }=\xi -%
\mathbb{E\xi }.%
\end{array}%
\right.  \label{defb}
\end{equation}%
Obviously, (\ref{defb}) is kind of decoupled FBSDEs whose solvability is
trivial.

\begin{lemma}
\label{applemma}Assume that \emph{(A1)} and \emph{(A2)} are in force, there
exists a positive constant $\delta _{0}\in \lbrack 0,1],$ such that if, a
priori, for some $\alpha _{0}\in \lbrack 0,1)$, for each $x_{0}\in R^{n},$ $%
\left( b_{0},\sigma _{0},\gamma _{0},\lambda _{0},\mu _{0},\psi _{0}\right)
\in L_{\mathbb{F}^{W}}^{2}(0,T;\mathbb{R}^{n}\times \mathbb{R}^{n}\times
\mathbb{R}^{n}\times \mathbb{R}^{n}\times \mathbb{R}^{n}\times \mathbb{R}%
^{n}),$ mean field FBSDEs (\ref{pfbsde}) have a unique adapted solution in $%
L_{\mathbb{F}^{W}}^{2}(0,T;\mathbb{R}^{n}\times \mathbb{R}^{n}\times \mathbb{%
R}^{n}\times \mathbb{R}^{n}\times \mathbb{R}^{n}\times \mathbb{R}^{n})$,
then for each $\delta \in \lbrack \alpha _{0},\alpha _{0}+\delta _{0}]$, for
each $x_{0}\in \mathbb{R}^{n},$ $\left( b_{0},\sigma _{0},\gamma
_{0},\lambda _{0},\mu _{0},\psi _{0}\right) \in L_{\mathbb{F}^{W}}^{2}(0,T;%
\mathbb{R}^{n}\times \mathbb{R}^{n}\times \mathbb{R}^{n}\times \mathbb{R}%
^{n}\times \mathbb{R}^{n}\times \mathbb{R}^{n}),$ Eq. (\ref{pfbsde}) also
have a unique solution in $L_{\mathbb{F}^{W}}^{2}(0,T;\mathbb{R}^{n}\times
\mathbb{R}^{n}\times \mathbb{R}^{n}\times \mathbb{R}^{n}\times \mathbb{R}%
^{n}\times \mathbb{R}^{n})$.
\end{lemma}

\paragraph{Proof.}

Define
\begin{equation*}
\mathcal{M}\left( 0,T\right) =L_{\mathbb{F}^{W}}^{2}(0,T;\mathbb{R}%
^{n}\times \mathbb{R}^{n}\times \mathbb{R}^{n})\times L_{\mathbb{F}^{W}}^{2,%
\mathcal{E}_{0}}(0,T;\mathbb{R}^{n}\times \mathbb{R}^{n})\times L_{\mathbb{F}%
^{W}}^{2}(0,T;\mathbb{R}^{n}).
\end{equation*}%
We set $\left( x^{0},\mathbb{E}x^{0},y^{0},\mathbb{E}%
y^{0},z^{0},p^{0},q^{0},k^{0}\right) =0,$ and solve iteratively the
following equations:%
\begin{equation}
\left\{
\begin{array}{l}
\mathrm{d}x^{i+1}=\big[\alpha \mathbf{B}\left( x^{i+1},\mathbb{E}%
x^{i+1},p^{i+1},q^{i+1},k^{i+1}\right) +\delta \mathbf{B}\left( x^{i},%
\mathbb{E}x^{i},p^{i},q^{i},k^{i}\right) +b_{0}\big ]\mathrm{d}t \\
\qquad \qquad +\left[ \alpha \Sigma \left( p^{i+1},q^{i+1},k^{i+1}\right)
+\delta \Sigma \left( p^{i},q^{i},k^{i}\right) +\sigma _{0}\right] \mathrm{d}%
W_{t}, \\
\mathrm{d}y^{i+1}=\big[\alpha \mathbf{F}\left( x^{i+1},\mathbb{E}%
x^{i+1},y^{i+1},\mathbb{E}y^{i+1},p^{i+1},q^{i+1},k^{i+1}\right) \\
\qquad \qquad +\delta \mathbf{F}\left( x^{i},\mathbb{E}x^{i},y^{i},\mathbb{E}%
y^{i},p^{i},q^{i},k^{i}\right) +\gamma -\mathbb{E\gamma }\big ]\mathrm{d}%
t+z^{i+1}\mathrm{d}W_{t}, \\
\mathrm{d}p^{i+1}=\big[\alpha \Xi \left( y^{i+1},\mathbb{E}%
y^{i+1},p^{i+1}\right) +\delta \Xi \left( y^{i},\mathbb{E}y^{i},p^{i}\right)
+\lambda _{0}-\mathbb{E}\lambda _{0}\big ]\mathrm{d}t \\
\mathrm{d}q^{i+1}=\big[\alpha \Upsilon \left( x^{i+1},\mathbb{E}%
x^{i+1},p^{i+1},q^{i+1}\right) +\delta \Upsilon \left( x^{i},\mathbb{E}%
x^{i},p^{i},q^{i}\right) \\
\qquad \qquad +\psi _{0}-\mathbb{E}\psi _{0}\big ]\mathrm{d}t+k^{i+1}\mathrm{%
d}W_{t}, \\
x_{0}=x\in \mathbb{R}^{n},\text{ }y_{T}^{i+1}=\alpha \Phi x_{T}^{i+1}+\delta
\Phi x_{T}^{i}+\upsilon -\mathbb{E\upsilon },\text{ }p_{0}=\varsigma -%
\mathbb{E\varsigma }, \\
q_{T}^{i+1}=\alpha \Phi ^{T}p_{T}^{i+1}-\alpha G\left( x_{T}^{i+1}-\mathbb{E}%
x_{T}^{i+1}\right) +\delta \Phi ^{T}\hat{p}_{T}^{i}-\delta G\left( x_{T}^{i}-%
\mathbb{E}x_{T}^{i}\right) +\xi -\mathbb{E\xi }.%
\end{array}%
\right.  \label{iter}
\end{equation}%
We set
\begin{eqnarray*}
\hat{x}^{i+1} &=&x^{i+1}-x^{i},\text{ }\hat{y}^{i+1}=y^{i+1}-y^{i}, \\
\hat{z}^{i+1} &=&z^{i+1}-z^{i},\text{ }\hat{p}^{i+1}=p^{i+1}-p^{i}, \\
\hat{q}^{i+1} &=&q^{i+1}-q^{i},\text{ }\hat{k}^{i+1}=k^{i+1}-k^{i}, \\
\mathbf{\hat{B}} &\mathbf{=}&\mathbf{B}\left( x^{i+1},\mathbb{E}%
x^{i+1},p^{i+1},q^{i+1},k^{i+1}\right) -\mathbf{B}\left( x^{i},\mathbb{E}%
x^{i},p^{i},q^{i},k^{i}\right) , \\
\hat{\Sigma} &=&\Sigma \left( p^{i+1},q^{i+1},k^{i+1}\right) -\Sigma \left(
p^{i},q^{i},k^{i}\right) , \\
\mathbf{\hat{F}} &=&\mathbf{F}\left( x^{i+1},\mathbb{E}x^{i+1},y^{i+1},%
\mathbb{E}y^{i+1},z^{i+1},p^{i+1},q^{i+1},k^{i+1}\right) \\
&&-\mathbf{F}\left( x^{i},\mathbb{E}x^{i},y^{i},\mathbb{E}%
y^{i},z^{i},p^{i},q^{i},k^{i}\right) , \\
\hat{\Xi} &=&\Xi \left( y^{i+1},\mathbb{E}y^{i+1}p^{i+1}\right) -\Xi \left(
y^{i},\mathbb{E}y^{i},p^{i}\right) , \\
\hat{\Upsilon} &=&\Upsilon \left( x^{i+1},\mathbb{E}x^{i+1},p^{i+1},q^{i+1}%
\right) -\Upsilon \left( x^{i},\mathbb{E}x^{i},p^{i},q^{i}\right)
\end{eqnarray*}%
with%
\begin{eqnarray*}
\hat{\varphi} &=&\varphi \left( p^{i+1},q^{i+1},k^{i+1}\right) -\varphi
\left( p^{i},q^{i},k^{i}\right) \\
&=&\mathbf{P}_{\Gamma }\left[ R_{t}^{-1}\left(
B_{t}^{T}q^{i+1}+K_{t}^{T}p^{i+1}+D_{t}^{T}k^{i+1}\right) \right] -\mathbf{P}%
_{\Gamma }\left[ R_{t}^{-1}\left(
B_{t}^{T}q^{i}+K_{t}^{T}p^{i}+D_{t}^{T}k^{i}\right) \right]
\end{eqnarray*}%
Now introduce a map $I_{\alpha _{0}}:\left(
x^{i},y^{i},z^{i},p^{i},q^{i},k^{i}\right) \rightarrow \left(
x^{i+1},y^{i+1},z^{i+1},p^{i+1},q^{i+1},k^{i+1}\right) \in \mathcal{M}\left(
0,T\right) $ by the following mean fields FBSDEs:%
\begin{equation}
\left\{
\begin{array}{l}
\mathrm{d}\hat{x}^{i+1}=\big[\alpha _{0}\mathbf{\hat{B}}\left( \hat{x}^{i+1},%
\mathbb{E}\hat{x}^{i+1},\hat{p}^{i+1},\hat{q}^{i+1},\hat{k}^{i+1}\right)
+\delta \mathbf{B}\left( \hat{x}^{i},\mathbb{E}\hat{x}^{i},\hat{p}^{i},\hat{q%
}^{i},\hat{k}^{i}\right) \big ]\mathrm{d}t \\
\qquad \qquad +\left[ \alpha _{0}\Sigma \left( \hat{p}^{i+1},\hat{q}^{i+1},%
\hat{k}^{i+1}\right) +\delta \Sigma \left( \hat{p}^{i},\hat{q}^{i},\hat{k}%
^{i}\right) \right] \mathrm{d}W_{t}, \\
\mathrm{d}\hat{y}^{i+1}=\big[\alpha _{0}\mathbf{\hat{F}}\left( \hat{x}^{i+1},%
\mathbb{E}\hat{x}^{i+1},\hat{y}^{i+1},\mathbb{E}\hat{y}^{i+1},\hat{z}^{i+1},%
\hat{p}^{i+1},\hat{q}^{i+1},\hat{k}^{i+1}\right) \\
\qquad \qquad +\delta \mathbf{\hat{F}}\left( \hat{x}^{i},\mathbb{E}\hat{x}%
^{i},\hat{y}^{i},\mathbb{E}\hat{y}^{i},\hat{z}^{i},\hat{p}^{i},\hat{q}^{i},%
\hat{k}^{i}\right) \big ]\mathrm{d}t+\hat{z}^{i+1}\mathrm{d}W_{t}, \\
\mathrm{d}\hat{p}^{i+1}=\big[\alpha _{0}\hat{\Xi}\left( \hat{y}^{i+1},%
\mathbb{E}\hat{y}^{i+1},\hat{p}^{i+1}\right) +\delta \hat{\Xi}\left( \hat{y}%
^{i},\mathbb{E}\hat{y}^{i},\hat{p}^{i}\right) \big ]\mathrm{d}t \\
\mathrm{d}\hat{q}^{i+1}=\big[\alpha _{0}\hat{\Upsilon}\left( \hat{x}^{i+1},%
\mathbb{E}\hat{x}^{i+1},\hat{p}^{i+1},\hat{q}^{i+1}\right) +\delta \hat{%
\Upsilon}\left( \hat{x}^{i},\mathbb{E}\hat{x}^{i},\hat{p}^{i},\hat{q}%
^{i}\right) \big ]\mathrm{d}t+\hat{k}^{i+1}\mathrm{d}W_{t}, \\
x_{0}=x\in \mathbb{R}^{n},\text{ }\hat{y}_{T}^{i+1}=\alpha \Phi \hat{x}%
_{T}^{i+1}+\delta \Phi \hat{x}_{T}^{i},\text{ }p_{0}=0, \\
\hat{q}_{T}^{i+1}=\alpha \Phi ^{T}\hat{p}_{T}^{i+1}-\alpha G\left( \hat{x}%
_{T}^{i+1}-\mathbb{E}\hat{x}_{T}^{i+1}\right) +\delta \Phi ^{T}\hat{p}%
_{T}^{i}-\delta G\left( \hat{x}_{T}^{i}-\mathbb{E}\hat{x}_{T}^{i}\right) .%
\end{array}%
\right.  \label{pri}
\end{equation}%
Applying the It\^{o}'s formula to $\left\langle \hat{x}^{i+1},\hat{q}%
^{i+1}\right\rangle -\left\langle \hat{y}^{i+1},\hat{p}^{i+1}\right\rangle $
on $\left[ 0,T\right] ,$ we have
\begin{eqnarray*}
&&\mathbb{E}\left[ \left\langle \hat{q}_{T}^{i+1},\hat{x}_{T}^{i+1}\right%
\rangle -\left\langle \hat{y}_{T}^{i+1},\hat{p}_{T}^{i+1}\right\rangle %
\right] \\
&=&\alpha _{0}\bigg [\mathbb{E}\int_{0}^{T}\left\langle \hat{q}^{i+1},%
\mathbf{\hat{B}}\left( \hat{x}^{i+1},\mathbb{E}\hat{x}^{i+1},\hat{p}^{i+1},%
\hat{q}^{i+1},\hat{k}^{i+1}\right) \right\rangle \\
&&+\left\langle \hat{x}^{i+1},\hat{\Upsilon}\left( \hat{x}^{i+1},\mathbb{E}%
\hat{x}^{i+1},\hat{p}^{i+1},\hat{q}^{i+1}\right) \right\rangle \\
&&+\left\langle \hat{k}^{i+1},\Sigma \left( \hat{p}^{i+1},\hat{q}^{i+1},\hat{%
k}^{i+1}\right) \right\rangle -\left\langle \hat{y}^{i+1},\hat{\Xi}\left(
\hat{y}^{i+1},\mathbb{E}\hat{y}^{i+1},\hat{p}^{i+1}\right) \right\rangle \\
&&-\left\langle \hat{p}^{i+1},\mathbf{\hat{F}}\left( \hat{x}^{i+1},\mathbb{E}%
\hat{x}^{i+1},\hat{y}^{i+1},\mathbb{E}\hat{y}^{i+1},\hat{p}^{i+1},\hat{q}%
^{i+1},\hat{k}^{i+1}\right) \right\rangle \mathrm{d}t\bigg ] \\
&&+\delta \bigg [\mathbb{E}\int_{0}^{T}\left\langle \hat{q}^{i+1},\mathbf{B}%
\left( \hat{x}^{i},\mathbb{E}\hat{x}^{i},\hat{p}^{i},\hat{q}^{i},\hat{k}%
^{i}\right) \right\rangle +\left\langle \hat{x}^{i+1},\hat{\Upsilon}\left(
\hat{x}^{i},\mathbb{E}\hat{x}^{i},\hat{p}^{i},\hat{q}^{i}\right)
\right\rangle \\
&&+\left\langle \hat{k}^{i+1},\Sigma \left( \hat{p}^{i},\hat{q}^{i},\hat{k}%
^{i}\right) \right\rangle -\left\langle \hat{y}^{i+1},\hat{\Xi}\left( \hat{y}%
^{i},\mathbb{E}\hat{y}^{i},\hat{p}^{i}\right) \right\rangle \\
&&-\hat{p}^{i+1},\mathbf{\hat{F}}\left( \hat{x}^{i},\mathbb{E}\hat{x}^{i},%
\hat{y}^{i},\mathbb{E}\hat{y}^{i},\hat{p}^{i},\hat{q}^{i},\hat{k}^{i}\right)
\mathrm{d}t\bigg ] \\
&=&\alpha _{0}\bigg [\mathbb{E}\int_{0}^{T}\left\langle \hat{p}^{i+1}+\hat{q}%
^{i+1}+\hat{k}^{i+1},\hat{\varphi}\left( \hat{p}^{i+1},\hat{q}^{i+1},\hat{k}%
^{i+1}\right) \right\rangle \\
&&+\left\langle \hat{x}^{i+1},Q\left( \hat{x}^{i+1}-\mathbb{E}\hat{x}%
^{i+1}\right) \right\rangle +\left\langle \hat{y}^{i+1},L\left( \hat{y}%
^{i+1}-\mathbb{E}\hat{y}^{i+1}\right) \right\rangle \mathrm{d}t\bigg ] \\
&&+\delta \bigg [\mathbb{E}\int_{0}^{T}\left\langle \hat{q}^{i+1},\delta
\mathbf{B}\left( \hat{x}^{i},\mathbb{E}\hat{x}^{i},\hat{p}^{i},\hat{q}^{i},%
\hat{k}^{i}\right) \right\rangle +\left\langle \hat{x}^{i+1},\hat{\Upsilon}%
\left( \hat{x}^{i},\mathbb{E}\hat{x}^{i},\hat{p}^{i},\hat{q}^{i}\right)
\right\rangle \\
&&+\left\langle \hat{k}^{i+1},\Sigma \left( \hat{p}^{i},\hat{q}^{i},\hat{k}%
^{i}\right) \right\rangle -\left\langle \hat{y}^{i+1},\hat{\Xi}\left( \hat{y}%
^{i},\mathbb{E}\hat{y}^{i},\hat{p}^{i}\right) \right\rangle \\
&&-\hat{p}^{i+1},\mathbf{\hat{F}}\left( \hat{x}^{i},\mathbb{E}\hat{x}^{i},%
\hat{y}^{i},\mathbb{E}\hat{y}^{i},\hat{p}^{i},\hat{q}^{i},\hat{k}^{i}\right)
\mathrm{d}t\bigg ].
\end{eqnarray*}

After simple computation, we have%
\begin{eqnarray*}
&&\mathbb{E}\left\langle \hat{x}_{T}^{i+1},\alpha _{0}G\left( \hat{x}^{i+1}-%
\mathbb{E}\hat{x}^{i+1}\right) \right\rangle +\alpha _{0}\bigg [\mathbb{E}%
\int_{0}^{T}\left\langle \hat{p}^{i+1}+\hat{q}^{i+1}+\hat{k}^{i+1},\hat{%
\varphi}\left( \hat{p}^{i+1},\hat{q}^{i+1},\hat{k}^{i+1}\right) \right\rangle
\\
&&+\left\langle \hat{x}^{i+1},Q\left( \hat{x}^{i+1}-\mathbb{E}\hat{x}%
^{i+1}\right) \right\rangle +\left\langle \hat{y}^{i+1},L\left( \hat{y}%
^{i+1}-\mathbb{E}\hat{y}^{i+1}\right) \right\rangle \mathrm{d}t\bigg ] \\
&=&-\delta \bigg [\mathbb{E}\int_{0}^{T}\left\langle \hat{q}^{i+1},\delta
\mathbf{B}\left( \hat{x}^{i},\mathbb{E}\hat{x}^{i},\hat{p}^{i},\hat{q}^{i},%
\hat{k}^{i}\right) \right\rangle +\left\langle \hat{x}^{i+1},\hat{\Upsilon}%
\left( \hat{x}^{i},\mathbb{E}\hat{x}^{i},\hat{p}^{i},\hat{q}^{i}\right)
\right\rangle \\
&&+\left\langle \hat{k}^{i+1},\Sigma \left( \hat{p}^{i},\hat{q}^{i},\hat{k}%
^{i}\right) \right\rangle -\left\langle \hat{y}^{i+1},\hat{\Xi}\left( \hat{y}%
^{i},\mathbb{E}\hat{y}^{i},\hat{p}^{i}\right) \right\rangle \\
&&-\hat{p}^{i+1},\mathbf{\hat{F}}\left( \hat{x}^{i},\mathbb{E}\hat{x}^{i},%
\hat{y}^{i},\mathbb{E}\hat{y}^{i},\hat{p}^{i},\hat{q}^{i},\hat{k}^{i}\right)
\mathrm{d}t\bigg ] \\
&&+\delta \left\langle \hat{x}_{T}^{i+1},\Phi ^{T}\hat{p}_{T}^{i}-G\left(
\hat{x}_{T}^{i}-E\hat{x}_{T}^{i}\right) \right\rangle -\delta \left\langle
\hat{x}_{T}^{i},\Phi ^{T}\hat{p}_{T}^{i+1}\right\rangle
\end{eqnarray*}%
By using the monotonicity property of $\varphi \left( p,q,k\right) $
(Proposition \ref{pro3} below and the classical geometric inequality and
Lipschitz property of projection operator (Proposition \ref{pro2}), it
follows that.
\begin{eqnarray}
&&\alpha _{0}\mathbb{E}\left\langle \hat{x}_{T}^{i+1},G\left( \hat{x}%
_{T}^{i+1}-\mathbb{E}\hat{x}_{T}^{i+1}\right) \right\rangle +\alpha _{0}%
\bigg [\mathbb{E}\int_{0}^{T}\left\langle \hat{x}^{i+1},Q\left( \hat{x}%
^{i+1}-\mathbb{E}\hat{x}^{i+1}\right) \right\rangle  \notag \\
&&+\left\langle \hat{y}^{i+1},L\left( \hat{y}^{i+1}-\mathbb{E}\hat{y}%
^{i+1}\right) \right\rangle \mathrm{d}t\bigg ]  \notag \\
&\leq &-\delta \bigg [\mathbb{E}\int_{0}^{T}\left\langle \hat{q}%
^{i+1},\delta \mathbf{B}\left( \hat{x}^{i},\mathbb{E}\hat{x}^{i},\hat{p}^{i},%
\hat{q}^{i},\hat{k}^{i}\right) \right\rangle +\left\langle \hat{x}^{i+1},%
\hat{\Upsilon}\left( \hat{x}^{i},\mathbb{E}\hat{x}^{i},\hat{p}^{i},\hat{q}%
^{i}\right) \right\rangle  \notag \\
&&+\left\langle \hat{k}^{i+1},\Sigma \left( \hat{p}^{i},\hat{q}^{i},\hat{k}%
^{i}\right) \right\rangle -\left\langle \hat{y}^{i+1},\hat{\Xi}\left( \hat{y}%
^{i},\mathbb{E}\hat{y}^{i},\hat{p}^{i}\right) \right\rangle  \notag \\
&&-\hat{p}^{i+1},\mathbf{\hat{F}}\left( \hat{x}^{i},\mathbb{E}\hat{x}^{i},%
\hat{y}^{i},\mathbb{E}\hat{y}^{i},\hat{p}^{i},\hat{q}^{i},\hat{k}^{i}\right)
\mathrm{d}t\bigg ]  \notag \\
&&+\delta \left\langle \hat{x}_{T}^{i+1},\Phi ^{T}\hat{p}_{T}^{i}-G\left(
\hat{x}_{T}^{i}-E\hat{x}_{T}^{i}\right) \right\rangle -\delta \left\langle
\hat{x}_{T}^{i},\Phi ^{T}\hat{p}_{T}^{i+1}\right\rangle  \notag \\
&\leq &\delta C_{1}\bigg [\mathbb{E}\int_{0}^{T}\left( \left\vert \hat{x}%
^{i}\right\vert ^{2}+\left\vert \hat{y}^{i}\right\vert ^{2}+\left\vert \hat{p%
}^{i}\right\vert ^{2}+\left\vert \hat{q}^{i}\right\vert ^{2}+\left\vert \hat{%
k}^{i}\right\vert ^{2}\right) \mathrm{d}t  \notag \\
&&+\delta C_{1}\bigg [\mathbb{E}\int_{0}^{T}\left( \left\vert \hat{x}%
^{i+1}\right\vert ^{2}+\left\vert \hat{y}^{i+1}\right\vert ^{2}+\left\vert
\hat{p}^{i+1}\right\vert ^{2}+\left\vert \hat{q}^{i+1}\right\vert
^{2}+\left\vert \hat{k}^{i+1}\right\vert ^{2}\right) \mathrm{d}t  \notag \\
&&+\delta C_{1}\mathbb{E}\left\vert \hat{x}_{T}^{i+1}\right\vert ^{2}+\delta
C_{1}\mathbb{E}\left\vert \hat{x}_{T}^{i}\right\vert ^{2}+\delta C_{1}%
\mathbb{E}\left\vert \hat{p}_{T}^{i+1}\right\vert ^{2}+\delta C_{1}\mathbb{E}%
\left\vert \hat{p}_{T}^{i}\right\vert ^{2}.  \label{ineq1}
\end{eqnarray}%
Next by baisc technique in SDE, BSDE, we have%
\begin{eqnarray}
&&\mathbb{E}\left[ \int_{0}^{T}\left\vert \hat{p}^{i+1}\right\vert ^{2}%
\mathrm{d}t\right] +\mathbb{E}\left\vert \hat{p}_{T}^{i+1}\right\vert ^{2}
\notag \\
&\leq &\delta C_{2}\mathbb{E}\int_{0}^{T}\left( \left\vert \hat{y}%
^{i}\right\vert ^{2}+\left\vert \hat{p}^{i}\right\vert ^{2}\right) \mathrm{d}%
s+C_{2}\mathbb{E}\int_{0}^{T}\left\vert L\left( \hat{y}^{i+1}-\mathbb{E}\hat{%
y}^{i+1}\right) \right\vert ^{2}\mathrm{d}s,  \label{ineq2}
\end{eqnarray}

\begin{eqnarray}
&&\mathbb{E}\left[ \int_{0}^{T}\left\vert \hat{q}^{i+1}\right\vert
^{2}dt+\int_{0}^{T}\left\vert \hat{k}^{i+1}\right\vert ^{2}\mathrm{d}s\right]
\notag \\
&\leq &\delta C_{3}\mathbb{E}\int_{0}^{T}\left( \left\vert \hat{x}%
^{i}\right\vert ^{2}+\left\vert \hat{p}^{i}\right\vert ^{2}+\left\vert \hat{q%
}^{i}\right\vert ^{2}\right) \mathrm{d}s+C_{3}\mathbb{E}\int_{0}^{T}\left%
\vert Q\left( \hat{x}^{i+1}-\mathbb{E}\hat{x}^{i+1}\right) \right\vert ^{2}%
\mathrm{d}s  \notag \\
&&+C_{3}\mathbb{E}\left\vert G\left( \hat{x}_{T}^{i+1}-\mathbb{E}\hat{x}%
_{T}^{i+1}\right) \right\vert ^{2}+\delta C_{3}\mathbb{E}\left[ \left\vert
\hat{p}_{T}^{i}\right\vert ^{2}+\left\vert \hat{x}_{T}^{i}\right\vert ^{2}%
\right] ,  \label{ineq3}
\end{eqnarray}%
and
\begin{eqnarray}
&&\mathbb{E}\int_{0}^{T}\left\vert \hat{x}_{t}^{i+1}\right\vert ^{2}\mathrm{d%
}t+\mathbb{E}\left\vert \hat{x}_{T}^{i+1}\right\vert ^{2}  \notag \\
&\leq &\delta C_{4}\mathbb{E}\int_{0}^{T}\left( \left\vert \hat{x}%
^{i}\right\vert ^{2}+\left\vert \hat{p}^{i}\right\vert ^{2}+\left\vert \hat{q%
}^{i}\right\vert ^{2}+\left\vert \hat{k}^{i}\right\vert ^{2}\right) \mathrm{d%
}s  \notag \\
&&+C_{4}\mathbb{E}\int_{0}^{T}\big (\left\vert \hat{p}^{i+1}\right\vert
^{2}+\left\vert \hat{q}^{i+1}\right\vert ^{2}+\left\vert \hat{k}%
^{i+1}\right\vert ^{2}\big )\mathrm{d}s.  \label{ineq4}
\end{eqnarray}%
Moreover,%
\begin{eqnarray}
&&\mathbb{E}\int_{0}^{T}\left( \left\vert \hat{y}^{i+1}\right\vert
^{2}+\left\vert \hat{z}^{i+1}\right\vert ^{2}\right) \mathrm{d}t  \notag \\
&\leq &\delta \mathbb{E}C_{5}\int_{0}^{T}\left( \left\vert \hat{x}%
^{i}\right\vert ^{2}+\left\vert \hat{y}^{i}\right\vert ^{2}+\left\vert \hat{p%
}^{i}\right\vert ^{2}+\left\vert \hat{q}^{i}\right\vert ^{2}+\left\vert \hat{%
k}^{i}\right\vert ^{2}\right) \mathrm{d}s  \notag \\
&&+C_{5}\mathbb{E}\int_{0}^{T}\big (\left\vert \hat{p}^{i+1}\right\vert
^{2}+\left\vert \hat{q}^{i+1}\right\vert ^{2}+\left\vert \hat{k}%
^{i+1}\right\vert ^{2}\big )\mathrm{d}s  \notag \\
&&+C_{5}\mathbb{E}\left\vert \hat{x}_{T}^{i+1}\right\vert ^{2}+\delta C_{4}%
\mathbb{E}\left\vert \hat{x}_{T}^{i}\right\vert ^{2}.  \label{ineq5}
\end{eqnarray}%
Observe that inequality (\ref{ineq2}) does not contain $\hat{x}^{i+1}$ and $%
\hat{y}^{i+1}.$ Combining (\ref{ineq1})-(\ref{ineq5}), by similar method
used in \cite{hp}, we have, for some $\delta _{0}\in \left( 0,1\right) ,$
\begin{eqnarray*}
&&\mathbb{E}\int_{0}^{T}\left( \left\vert \hat{p}^{i+1}\right\vert
^{2}+\left\vert \hat{q}^{i+1}\right\vert ^{2}+\left\vert \hat{k}%
^{i+1}\right\vert ^{2}\right) \mathrm{d}s+\mathbb{E}\left\vert \hat{p}%
_{T}^{i+1}\right\vert ^{2} \\
&\leq &\delta _{0}C_{5}\mathbb{E}\int_{0}^{T}\left( \left\vert \hat{p}%
^{i}\right\vert ^{2}+\left\vert \hat{q}^{i}\right\vert ^{2}+\left\vert \hat{k%
}^{i}\right\vert ^{2}\right) \mathrm{d}s+\mathbb{E}\left\vert \hat{p}%
_{T}^{i}\right\vert ^{2},
\end{eqnarray*}%
which means that the map $I_{\alpha _{0}+\delta _{0}}$ is a contraction.
\hfill $\Box $

(\textbf{Existence}) We can solve Eq. (\ref{pfbsde}) successively for the
case $\alpha \in \left[ 0,\delta _{0}\right] ,$ $\left[ \delta _{0},2\delta
_{0}\right] ,\cdots $When $\alpha =1,$ we deduce immediately that the
solution to Eq. (\ref{cc}) exists. \hfill $\Box $

\section{Properties of projection}

\label{app2}

We recall the following properties of projection $\mathbf{P}_{U}$ onto a
closed convex set $U$, see \cite{Brezis 2010}, Chapter 5.

\begin{theorem}
\label{projection theorem} For a nonempty closed convex set $U\subset
\mathbb{R}^{m}$, for every $x\in \mathbb{R}^{m}$, there exists a unique $%
x^{\ast }\in U$, such that
\begin{equation*}
|x-x^{\ast }|=\min_{y\in \Gamma }|x-y|=:dist(x,U).
\end{equation*}%
Moreover, $x^{\ast }$ is characterized by the property
\begin{equation}
x^{\ast }\in U,\quad \big<x^{\ast }-x,x^{\ast }-y\big>\leq 0\qquad \forall
y\in U.  \label{projection characterization}
\end{equation}%
The above element $x^{\ast }$ is called the projection of $x$ onto $U$ and
is denoted by $\mathbf{P}_{U}[x]$.
\end{theorem}

One can immediately obtain the following

\begin{proposition}
\label{pro1}Let $U\subset \mathbb{R}^{m}$ be a nonempty closed convex set,
then we have
\begin{equation}
\big|\mathbf{P}_{U}[x]-\mathbf{P}_{U}[y]\big|^{2}\leq \big<\mathbf{P}_{U}[x]-%
\mathbf{P}_{U}[y],x-y\big>.  \label{projection inequality}
\end{equation}
\end{proposition}


\begin{proposition}
\label{pro2}Let $U\subset \mathbb{R}^{m}$ be a nonempty closed convex set,
then the projection $\mathbf{P}_{U}$ does not increase the distance, i.e.
\begin{equation}
\big|\mathbf{P}_{U}[x]-\mathbf{P}_{U}[y]\big|\leq \big|x-y\big|.
\label{projection
lipschitz}
\end{equation}
\end{proposition}


\noindent Now let us consider $\mathbb{R}^{m}$ and the projection $\mathbf{P}%
_{U}$ both with the norm $\Vert \cdot \Vert _{R_{0}}:=\langle R_{0}^{\frac{1%
}{2}}\cdot ,R_{0}^{\frac{1}{2}}\cdot \rangle $, from (\ref{projection
inequality}), we have

\begin{proposition}
\label{pro3}Let $U\subset \mathbb{R}^{m}$ be a nonempty closed convex set,
then
\begin{equation*}
\langle \langle \mathbf{P}_{U}[x]-\mathbf{P}_{U}[y],x-y\rangle \rangle
=\left\langle R^{\frac{1}{2}}\bigg(\mathbf{P}_{U}[x]-\mathbf{P}_{U}[y]\bigg)%
,R^{\frac{1}{2}}(x-y)\right\rangle \geq 0.
\end{equation*}
\end{proposition}

\noindent The proofs of Proposition \ref{pro1}-Proposition \ref{pro3} can be
found in \cite{Ba, Brezis 2010}.

\noindent \textbf{Conflict of Interest:} The authors declare that they have
no conflict of interest.

\noindent \textbf{Acknowledgements.} The authors wish to thank the editors
and two referees for their valuable comments and constructive suggestions
which improved the presentation of this manuscript.

\end{document}